\documentclass[11pt, reqno]{amsart}

\usepackage{amsmath, amsthm, amssymb, bbm, graphicx, hyperref, epsfig, amsfonts, bm}
\usepackage[usenames,dvipsnames]{xcolor}

\usepackage{tikz-cd}
\usepackage{mathtools}

\usepackage{graphicx}
\usepackage{listings}

\usepackage[normalem]{ulem}

\usepackage[margin=1in]{geometry}
\usepackage{lstautogobble}
\usepackage{enumerate}
\usepackage[shortlabels]{enumitem}
\usepackage{thmtools}
\usepackage{thm-restate}
\usepackage{amsthm}
\usepackage{verbatim}
\usepackage{accents}
\usepackage{mathtools}
\usepackage{physics}

\usepackage[capitalise]{cleveref}  

\usepackage[bottom]{footmisc}
\crefformat{equation}{(#2#1#3)}

\usepackage{float}
\restylefloat{table}

\usepackage{mathrsfs}
\setlist{  
  listparindent=\parindent,
  parsep=0pt,
}

\theoremstyle{plain}

\newtheorem{thm}{Theorem}
\newtheorem{prop}{Proposition}[section]
\newtheorem{lemma}[prop]{Lemma}
\newtheorem{cor}[prop]{Corollary}

\theoremstyle{definition}
\newtheorem{defn}[prop]{Definition}
\newtheorem{remark}[defn]{Remark}

\numberwithin{equation}{section}

\Crefname{thm}{Theorem}{Theorems}
\Crefname{prop}{Proposition}{Propositions}

\numberwithin{equation}{section} 

\DeclareMathOperator*{\argmin}{arg\,min}

\newcommand{\eps}{\varepsilon}

\newcommand{\diam}{\operatorname{diam}}
\newcommand{\Var}{\mathrm{Var}}

\newcommand{\be}{\begin{equation}}
\newcommand{\ee}{\end{equation}}
\newcommand{\bea}{\begin{eqnarray}}
\newcommand{\eea}{\end{eqnarray}}
\newcommand{\beas}{\begin{eqnarray*}}
\newcommand{\eeas}{\end{eqnarray*}}

\newcommand{\inner}[2]{\langle #1, \ #2 \rangle}

\newcommand{\co}{\mathrm{c}}

\newcommand{\TV}{\operatorname{TV}}

\newcommand{\bR}{\mathbb{R}}
\newcommand{\bN}{\mathbb{N}}
\newcommand{\bP}{\mathbb{P}}
\newcommand{\bE}{\mathbb{E}}

\newcommand{\bdel}{\mathbf{\Delta}}

\newcommand{\cG}{\mathcal{G}}

\newcommand{\cP}{\mathcal{P}_{q}}
\newcommand{\cPm}{\mathcal{P}_{q}^{\otimes m}}

\newcommand{\fT}{\mathbf{T_2}}
\newcommand{\fS}{\mathbf{S}}
\newcommand{\fA}{\mathbf{A}}
\newcommand{\fB}{\mathbf{B}}
\newcommand{\fC}{\mathbf{C}}
\newcommand{\fI}{\mathbf{I}}
\newcommand{\fe}{\mathbf{e}}
\newcommand{\fJ}{\mathbf{J}}

\newcommand{\fx}{\mathbf{x}}
\newcommand{\fz}{\mathbf{z}}
\newcommand{\fy}{\mathbf{y}}
\newcommand{\fw}{\mathbf{w}}

\newcommand{\crc}{\operatorname{cr}}

\newcommand{\bp}{\begin{proof}}
\newcommand{\ep}{\end{proof}}

\newcommand{\mix}{\operatorname{MIX}}

\newcommand{\CWP}{\operatorname{CWP}}

\begin{document}
\nocite{glauberising, completepotts, Glauber, yang, twocomp, fluctuising, fluctupotts, LDPblock, limitbipartite, ellislimit, markov, threegeneral, blumepath, rapid, aggregate, landscapepotts, heejune, path, dingising, ellisentropy}
\title[Dynamical Phase Transition for the homogeneous multi-component Potts model]{Dynamical Phase Transition for the homogeneous multi-component Curie-Weiss-Potts model}

\author{Kyunghoo Mun}
\address{Department of Mathematical Sciences, Carnegie Mellon University, Pittsburgh, PA 15213}
\email{kmun@andrew.cmu.edu}


\begin{abstract}
In this paper, we study the homogeneous multi-component Curie-Weiss-Potts model with $q \geq 3$ spins. The model is defined on the complete graph $K_{Nm}$, whose vertex set is equally partitioned into $m$ components of size $N$. For a configuration $\sigma: \{1, \cdots, Nm\} \to \{1, \cdots, q\},$ the Gibbs measure is defined by \begin{align*}
    \mu_{N, \beta}(\sigma) = {1 \over Z_{N, \beta}} \exp( {\beta \over N} \sum_{v, w =1}^{Nm} \mathcal{J}(v, w) \mathbbm{1}\{ \sigma(v) = \sigma(w)\}),
\end{align*} where $Z_{N, \beta}$ is a normalizing constant, and $\beta>0$ is the inverse temperature parameter. The interaction coefficient is \begin{align*}
    \mathcal{J}(v, w) = J \cdot  \begin{cases}
        \frac{1}{1+(m-1)\lambda} & \text{if } v, w \text{ are in the same component,}\\
        \frac{\lambda}{1+(m-1)\lambda} & \text{if } v, w \text{ are in different components,}
    \end{cases} 
\end{align*} with $\lambda \in (0, 1),$ the relative strength of inter-component interaction to intra-component interaction, and $J>0,$ the effective interaction strength. We identify a dynamical phase transition at the critical inverse temperature $\beta_{\crc} = \beta_{s}(q)/J$, where $\beta_{s}(q)$ is maximal inverse temperature guaranteeing a unique critical point of the free energy in the Curie-Weiss-Potts model \cite{glauberising}. By extending the aggregate path method \cite{rapid}  to our multi-component setting, we prove $O(N \log N)$ mixing time in the high-temperature regime $\beta<\beta_{s}(q)/J.$ In the low-temperature regime $\beta > \beta_{s}(q)/J,$ we further show exponential mixing time by a metastability. This is the first result for the dynamical phase transition in the multi-component Potts model.
\end{abstract}

\date{\today}

\maketitle

\section{Introduction} \label{intro} 
\subsection{Background and proof idea}
Dynamical phase transition has been a key topic in understanding how the convergence rate of a Markov process depends on model parameters such as the inverse temperature \(\beta>0\). Typically, the Ising model (the fundamental model for magnetization in statistical physics) describes the evolution of configurations of two-spin. In particular, the {Curie-Weiss-Ising} model, an Ising model on the complete graph \(K_{N}\) with \(N\) vertices, has been extensively studied in the mathematical literature \cite{dingising, glauberising, markov}. When temperature is high ($\beta\ll1$), entropy dominates the interaction energy and the system mixes rapidly. In contrast, in the low‐temperature regime ($ \beta\gg1$), an energy barrier prevents the a system from mixing fastly, resulting in slow mixing. This heuristic mixing change can be mathematically formulated as a dynamical phase transition, characterized by critical inverse temperature $\beta_{\crc},$ at which the order of mixing time in terms of the number of vertices changes drastically.

A dynamical phase transition in the Curie-Weiss-Ising model with $\beta_{\crc}=1$ was studied in \cite{glauberising}.  They showed that the mixing time varies by $\beta>0$ as
\begin{enumerate}\label{DPP Ising}
    \item In the high‐temperature regime $\beta<1$,
    \begin{align} \label{DPP Ising 1}
        t_{\mix} = O(N\log N).
    \end{align}
    \item At the critical temperature $\beta=1$,
    \begin{align} \label{DPP Ising 2}
        t_{\mix} = \Theta(N^{3/2}).
    \end{align}
    \item In the low‐temperature regime $\beta>1$,
    \begin{align} \label{DPP Ising 3}
        t_{\mix} = \exp\bigl(\Omega(N)\bigr).
    \end{align}
\end{enumerate}
In the high-temperature regime, they also proved a cutoff phenomenon -- an abrupt convergence of a system to its stationary measure -- with
$t_{\mix}
    = [2(1-\beta)]^{-1} N \log N \pm O(N).$ Since there are only two spins in Ising model, the proportion of one spin uniquely determines the proportion of the other.  This allowed \cite{glauberising} to construct a monotone coupling, that preserves the order between proportions of two spins. Using this coupling, they obtained the $O(N\log N)$ upper bound when $\beta<1$.

{Curie-Weiss-Potts (CWP) model} generalizes the Curie-Weiss-Ising model by extending the number of spins from two to $q\ge3$. This extension significantly complicates the free energy landscape. While the phase transition is continuous when $q=2$, \cite{ellislimit} showed it becomes discontinuous for $q\ge3$. In fact, for $q\ge3,$ the critical point threshold
\begin{align} \label{beta_s def}
\beta_s(q)
\nonumber&=\sup\{\beta>0:\text{the free energy has a unique critical point}\}\\
&=\sup{\left\{\beta >0: \big(1+(q-1)e^{2\beta{1-qx \over q-1}}\big)^{-1}-x \neq 0 \ \textnormal{for all} \ x\in (1/q, 1) \right\}}
\end{align}
and the global minimizer threshold
\begin{align} \label{beta_c def}
\beta_c(q)
\nonumber&=\sup\{\beta>0:\text{the free energy has a unique global minimizer}\}\\
&=\frac{q-1}{q-2} \log(q-1)
\end{align}
satisfy $\beta_s(q)<\beta_c(q),$ while $\beta_s(2)=\beta_c(2)$.  Moreover, \cite{Glauber}  revealed that a dynamical phase transition occurs in the Curie-Weiss-Potts model at $\beta_{\crc}=\beta_s(q).$
\begin{enumerate}
  \item In the high-temperature regime $\beta<\beta_s(q)$, $$t_{\mix}=O(N\log N).$$
  \item At the critical temperature $\beta=\beta_s(q)$, $$t_{\mix}=\Theta(N^{4/3}).$$
  \item In the low‐temperature regime $\beta>\beta_s(q)$, $$t_{\mix}=\exp(\Omega(N)).$$
\end{enumerate}
They also proved a cutoff for $\beta<\beta_s(q)$ with
$
t_{\mix}
=[2(1-2\beta/q)]^{-1} N \log N \pm O(N).$ Since in the Potts model, the full proportion vector is not determined by a single spin proportion, they introduced more sophisticated couplings to get the $O(N\log N)$ upper bound.

Studies on Curie-Weiss type models have recently been extended to multi‐component settings by allowing more general interactions via component. Suppose the complete graph $K_{N}$ on vertex set $N$ is partitioned into $m$ components of fixed size (independent of $N$). The interaction coefficient between two vertices depends only on the components to which they belong. For the general multi-component Curie-Weiss-Ising model, the dynamical phase transition was verified in \cite{yang}. The results mirror to \eqref{DPP Ising 1}-\eqref{DPP Ising 3} in the Curie-Weiss-Ising model, upon replacing $\beta$ by $\beta \rho,$ where $\rho$ is determined by ratios of component sizes and interaction coefficients. (See \cite{yang} for details.) They generalized the monotone coupling  \cite{glauberising} to the multi-component Ising setting. However, the dynamical phase transition for the multi-component Curie-Weiss-Potts model has not yet been identified.  

For simpler analysis, homogeneity of component interactions was assumed in \cite{LDPblock, fluctupotts, twocomp}. Each components are supposed to have the same number of vertices. All intra-component interactions (between vertices in the same component) are equal, and all inter-component interactions (between vertices in different components) are equal, regardless of which components are chosen. In this homogeneous setting, \cite{LDPblock} used large deviation theory to identify the phase transition, and \cite{fluctupotts} proved a central limit theorem for sufficiently small $\beta$. Moreover, \cite{twocomp} studied the case $m=2$, $q=3$, $J=1,$ showing that the critical-point threshold $$ \sup \{ \beta >0: \text{the free energy of this model has a unique critical point}\} $$ coincides with $\beta_s(3)$ for the original Curie-Weiss-Potts model.

\cite{rapid} introduced the aggregate path method for the Curie–Weiss–Potts model, obtaining an upper bound $O(N\log N)$ when $\beta<\beta_{s}(q).$
They constructed small monotone segments and aggregated them, ensuring that the resulting coupling path is close to the original configuration path. \cite{aggregate} extended this method to the bipartite Curie-Weiss-Potts model ($m=2$ with zero intra-component interactions).

At low temperatures, local minimizers of the free energy emerge, a phenomenon known as metastability. \cite{landscapepotts} analyzed how metastability develops as $\beta$ increases, providing a complete description of all local minimizers in the Curie–Weiss–Potts model. In the homogeneous multi-component setting with $m=2$ and $q=3$, \cite{twocomp} further characterizes the free‐energy landscape, stating that the critical values of $\beta$ -- at which the number of critical points changes -- are independent of the relative strength $\lambda,$ defined as a ratio of inter‐component to intra‐component interaction.

\medskip

\subsection{Homogeneous multi-component Potts model}We study a homogeneous multi-component Curie-Weiss-Potts model on the complete graph $K_{Nm}=(V, E),$ where $V$ denotes the vertex set and $E$ is the edge set. We partition $V$ into $m$ identical components of size $N$ as \begin{align*}
    V = \bigcup_{i=1}^{m} V_{i}, \quad V_{i} = \{1, 2, \cdots, N\} \quad \forall 1 \leq i \leq m.
\end{align*} We emphasize that the number of components $m$ is independent of $N$. \footnote{An alternative convention lets the total number of vertices be $N$. To avoid the technicality that $N/m$ may not be integer, we instead let $N$ denote the size of each component.} A configuration $\sigma$ maps each vertex $v \in V$ to a spin $\sigma(v) \in \{1, 2, \cdots, q\}.$ We denote the collection of all configurations by $\Sigma = \{1, 2, \cdots, q\}^{V}.$ We write $\sigma=(\sigma^{1}, \sigma^{2}, \cdots, \sigma^{m}),$ where $\sigma^{i} = \sigma|_{V_i}$ is the restriction of $\sigma$ on  component $V_i$. The Hamiltonian $H_{N}: \Sigma \to \bR$ is
\begin{align*}
H_{N}(\sigma)&= =-{1 \over N} \sum_{v, w =1}^{Nm} \mathcal{J}(v, w) \mathbbm{1} \{ \sigma(v) = \sigma(w)\},
\end{align*}
where interaction $\mathcal{J} \in \bR^{Nm \times Nm}$ has the block structure \begin{align*}
    \mathcal{J} = \begin{bmatrix}
        &\fJ^{11} \mathbf{1}_{N} \ &\fJ^{12} \mathbf{1}_{N} \  &\cdots \ &\fJ^{1m} \mathbf{1}_{N} \\
        &\fJ^{21} \mathbf{1}_{N} \ &\fJ^{22} \mathbf{1}_{N} \  &\cdots \ &\fJ^{2m} \mathbf{1}_{N} \\ &\vdots \ &\vdots \ &\ddots \ &\vdots  \\
        &\fJ^{m1} \mathbf{1}_{N} \ &\fJ^{m2} \mathbf{1}_{N} \  &\cdots \ &\fJ^{mm} \mathbf{1}_{N}
    \end{bmatrix}.
\end{align*} Here  $\mathbf{1}_{N}$ is $N$ by $N$ matrix of all ones and $\mathbf{J} = (\fJ^{i_1 i_2}) \in \mathbb{R}^{m \times m}$ expresses the interactions among components. For simplicity, we assume homogeneity on $\mathcal{J}$ by setting intra-component and inter-component interactions as  \begin{align*}
    \fJ^{i i} = \frac{J}{1+(m-1) \lambda}, \text{ and } \quad \fJ^{i_1 i_2} = \frac{J \lambda}{1+(m-1) \lambda}, \quad \text{for } i_{1} \neq i_{2}, \ i \in [m].
\end{align*} Equivalently, \begin{align*}
    \mathbf{J}= J \bm{\Lambda} = \frac{J}{1 + (m-1)\lambda}
    \begin{bmatrix}
1 \ \; \ \lambda\; \ \cdots \ \; \ \lambda\\
\lambda \ \; \ 1 \; \ \cdots \ \; \ \lambda\\
\ddots\\
\lambda\ \; \ \lambda \; \ \cdots \ \; \ 1
\end{bmatrix},
\end{align*}
 where $J>0$ is the effective interaction coefficient, and $\lambda>0$ is the relative strength of inter-component interaction compared to intra-component interaction. 

 \begin{remark}
     In earlier studies, the effective interaction $J$ was set to $1.$ The original Curie-Weiss-Potts model \cite{Glauber, rapid} corresponds to $m=1.$ The bipartite Curie-Weiss-Potts model \cite{aggregate, limitbipartite} arises when $m=2$ and $\lambda \to \infty$. In \cite{twocomp}, the case $m=2,$ $q=3$ with $\lambda \in (0, 1)$ was analyzed. To ensure positive-definiteness of $\fJ,$ \cite{LDPblock, fluctupotts} assumed $\lambda \in (0, 1).$  
 \end{remark}

 From now on, we assume the interaction matrix $\fJ$ is positive-definite, i.e. $\lambda <1$, so that intra-component interaction is stronger than inter-component interaction. The Gibbs measure $\mu_{N, \beta}$ at inverse temperature $\beta>0$ is defined as 
\begin{align*}
\mu_{N, \beta}(\sigma)&={1 \over Z_{N, \beta}} \exp(-\beta H_{N}(\sigma))\\ 
&= {1 \over Z_{N, \beta}} \exp({\beta \over N} \sum_{i_1, i_2 = 1}^{m} \sum_{v, w=1}^{N} \mathbf{J}^{i_1 i_2} \mathbbm{1}\{\sigma^{i_1}(v)=\sigma^{i_2}(w)\}),
\end{align*}
where $Z_{N, \beta}=\sum_{\sigma \in \Sigma} \exp\big(- \beta H_{N}(\sigma)\big)$ is the partition function.

\medskip

\subsection{Glauber dynamics and mixing time}
We run the discrete time Glauber dynamics on the configuration space $\Sigma$. This defines a Markov chain $(\sigma_t)_{t \geq 0}$ with initial state $\sigma_0.$ At each time $t,$ given the current configuration $\sigma =\sigma_t$, we select a vertex $v \in V$ uniformly at random, and update its spin by sampling from the conditional Gibbs measure \begin{align*}
    \mu_{N, \beta}\big[ \sigma_{t+1}  \big|   \sigma_{t+1} (w) =  \sigma_{t}(w), \ \forall w \neq v \big].
\end{align*}

We denote $\bE_{\sigma}$ and $\bP_{\sigma},$ the expectation and probability measure when the chain starts from $\sigma_0 \in \Sigma$. It is well known that $(\sigma_t)_{t \geq 0}$ is reversible Markov process with the stationary measure $\mu_{N, \beta}.$ (See \cite{markov} for details.)

A mixing time is the smallest time that total-variation distance between distribution of $\sigma_t$ and the Gibbs measure is at most $\eps.$ Formally, for $\eps \in (0, 1),$ \begin{align*}
    t_{\mix}(\eps)= \inf \big\{ t>0: \max_{\sigma \in \Sigma} \| \bP_{\sigma}[\sigma_t \in \cdot] - \mu_{N, \beta} \|_{\TV} \leq \eps \big\}.
\end{align*} We set $t_{\mix} = t_{\mix}(1/4)$ and our main interest is its asymptotic behavior as $N \to \infty.$

\medskip

\subsection{Main Results} We show that the homogeneous multi-component Curie-Weiss-Potts model exhibits a dynamical phase transition at the critical inverse temperature \begin{align*}
    \beta_{\crc} = {\beta_{s}(q) \over J}.
\end{align*} Regardless of the relative strength $\lambda \in (0, 1),$ only the effective coefficient $J$ and the number of spins $q$ determine $\beta_{\crc}.$

In the high-temperature regime $\beta<\beta_{s}(q)/J,$ we prove a fast mixing. 

\begin{thm} \label{thrm 1}
In the homogeneous multi-component Curie-Weiss-Potts model, if $\beta<\beta_{s}(q)/J$, then for big enough $N>0,$ \begin{align*}
t_{\mix}= O(N \log N).
\end{align*} 
\end{thm}

In the low-temperature regime $\beta>\beta_{s}(q)/J,$ at least $q$ local minimizers of the free energy emerge. Escaping from one local minimizer to another then requires exponential time, and consequently the mixing time grows exponentially.   

\begin{thm} \label{thrm 2}
In the homogeneous multi-component Curie-Weiss-Potts model, for $\beta>\beta_{s}(q)/J$, we have exponential (slow) mixing as\begin{align*}
t_{\mix}= \exp( \Omega(N)).
\end{align*} 
\end{thm}
\begin{remark}
    At the critical temperature $\beta = \beta_{s}(q)/J,$ the order of mixing time remains an open problem. According to \cite{Glauber} and the Ising case \cite{glauberising, yang}, we conjecture that the mixing is of order $\Theta(N^{4/3}),$ matching the order in the Curie-Weiss-Potts model.
\end{remark}

\medskip

\subsection{Outline of the proof}We extend the aggregate path method to the homogeneous multi‐component setting. Our generalization is the first result of $O(N  \log N)$ mixing time for non-trivial component interactions (both inter-component and intra-component interactions are not assumed to be zero). We couple two copies of Glauber chain by selecting the same vertex $v \in V$ at each step and updating its spin independently. We show that this greedy coupling contracts the two chains, if the proportion of one initial configuration was close to the uniform distribution of spins. We approximate evolution of proportions  by connecting small monotone segments. On each segment, a single coordinate of the proportion matrix moves monotonically toward. Under the assumption $\beta < \beta_{s}(q)/J,$ we show that the proportion evolution gets closer to the uniform proportion in the $L^{1}$ norm. Combining this with the large deviation result, we obtain $O(N \log N)$ mixing time when $\beta< \beta_{s}(q)/J.$

In the Curie–Weiss–Potts model, when $\beta>\beta_{s}(q),$ a one-spin dominant configuration occurs as a local minimizer of the free energy. By concatenating this proportion vector identically across each component, we construct a local minimizer of the free energy in our multi-component model. We permute spin labels yielding $q$ local minimizers when $\beta>\beta_{s}(q)/J$. Applying Cheeger’s inequality (see \cite[Chapter 7]{markov}, \cite{Glauber}) we obtain the exponential lower bound
$\exp(\Omega(N))$ when $\beta >\beta_{s}(q)/J.$

\subsection{Organization of the paper}
This paper is organized as follows. In Section ~\ref{intro}, we introduce our model, define the Glauber dynamics and state our results. In Section \ref{pre}, we analyze the update probabilities for a configuration and present the preliminary large deviation results. In Section \ref{high}, we construct a greedy coupling and develop the   aggregate path method resulting fast mixing in the high-temperature regime. Finally, in Section \ref{low}, we prove metastability in the low-temperature regime $\beta>\beta_{s}(q)/J$.

\subsection{Acknowledgements}
I would sincerely appreciate Insuk Seo for introducing me to interacting particle systems and for many insightful discussions. This work was partially supported by the Seoul National University student-directed undergraduate research program (2022).

\section{Preliminaries} \label{pre} 
\subsection{Notations}In this paper, we use the following notations.

\begin{itemize}
    \item $f(N) = O(g(N)),$ if there exists a constant $c>0$ such that $f(N) \leq c g(N).$
    \item $f(N) = \Omega(g(N)),$ if there exists a constant $c>0$ such that $f(N) \geq c g(N).$
    \item $f(N) = \Theta(g(N)),$ if there exists a constant $c, C>0$ such that $cg(N) \leq f(N) \leq Cg(N).$
    \item For $n \in \bN,$ $[n]= \{1,  2,  \cdots,  n \}.$
    \item Given a matrix $\mathbf{A}$, $\mathbf{A}^{\intercal}$ denotes its transpose, $\mathbf{A}^{\co}$ denotes its complement, and $\operatorname{Tr}(\mathbf{A})$ denotes its trace.
    \item For matrix $\mathbf{B} \in \mathbb{R}^{m \times n}$ and $i \in [m], \ j \in [n],$ the entry located on the position $(i, j)$ is $\mathbf{B}^{ij}.$ We write \begin{align*}
    \mathbf{B}^{i*}=[\mathbf{B}^{i1}, \mathbf{B}^{i2}, \cdots, \mathbf{B}^{in}], \quad \mathbf{B}^{*j}=[\mathbf{B}^{1j}, \mathbf{B}^{2j}, \cdots, \mathbf{B}^{mj}]^{\intercal}.
    \end{align*}
    \item For $\mathbf{C} \in \mathbb{R}^{m \times n}$, the entrywise norms are \begin{align*}
    \|\mathbf{C}\|_{1}=\sum_{i=1}^{m} \sum_{j=1}^{n} |\mathbf{C}^{ij}|, \quad
    \|\mathbf{C}\|_{2} = \Big[ \sum_{i=1}^{m} \sum_{j=1}^{n} \big(\mathbf{C}^{ij}\big)^{2} \Big]^{1/2}. 
    \end{align*} The operator norms are
    \begin{align*}
    \| \fC\|_{1 \to 1} = \sup \{ \|\fC x\|_{1}: x \in \bR^{n},  \|x\|_{1} \leq 1 \} = \max_{1 \leq j \leq n} \sum_{i=1}^{m} |\fC^{ij}|,
    \end{align*}
    and $\|\mathbf{C}\|_{2 \to 2}$ is the spectral norm, which is the largest singular value of matrix $\fC.$
    \item We define an inner product between two matrices $\fA, \fB \in \bR^{m \times n}$ by \begin{align*}
        \inner{\fA}{\fB} = \sum_{i=1}^{m} \sum_{j=1}^{n} \fA^{ij} \fB^{ij}
    \end{align*} 
\end{itemize}
If it is not stated otherwise, we identify $\bR^{m \times n} \simeq \bR^{mn}$ and apply vector operations as above.

\medskip

\subsection{Proportion matrix and free energy in the homogenous multi-component Potts model}The configuration space has cardinality $|\Sigma| = q^{Nm}.$ We reduce a complexity by mapping each configuration $\sigma \in \Sigma$ to its \textit{proportion matrix} $\fS(\sigma),$ where for each $i \in [m] \ \text{and } j \in [q]$, \begin{align*}
     \fS(\sigma)^{ij} = {1 \over N} \sum_{v \in V_i} \mathbbm{1}_{\{\sigma(v)=j\}}.
\end{align*} We can view $\fS$ as the operator $\fS: \Sigma \to \mathcal{P}_{q}^{\otimes m},$ where \begin{align*}
\mathcal{P}_{q}^{\otimes m}= \Bigg\{\mathbf{z} \in \mathbb{R}^{m \times q}: \ {\mathbf{z}}^{ij} \geq 0, \quad \sum_{j=1}^{q}\mathbf{z}^{ij}=1 \ \text{for all } i \in [m] \Bigg\}.
\end{align*}
In these terms, the Gibbs measure $\mu_{N, \beta}(\sigma)$ becomes \begin{align} \label{Gibbs expression}
    \mu_{N, \beta}(\sigma)={1 \over Z_{N, \beta}} \exp \big({\beta N} \operatorname{Tr}((\fS(\sigma))^{\intercal}\fJ \fS(\sigma)) \big).
\end{align}
This induces a measure $\nu_{N, \beta}$ in $\cPm$ via \begin{align*}
    \nu_{N, \beta}(\fz) = \sum_{\sigma \in \Sigma:  \fS(\sigma) = \fz} \mu_{N,\beta} (\sigma) = {1 \over Z_{N, \beta}} \sum_{\sigma \in \Sigma:  \fS(\sigma) = \fz}  \exp \big({\beta N} \operatorname{Tr}(\fz^{\intercal}\fJ \fz) \big).
\end{align*}
The following lemma expresses $\nu_{N, \beta}$ by the free energy $F_{\beta}$ and a remainder $F^{\textnormal{sub}}.$

\begin{lemma}
    For every $\fz \in \cPm,$ \begin{align*}
        \nu_{N, \beta}(\fz) = {1 \over \Hat{Z}_{N, \beta}} \exp( -N \Big(F_{\beta}(\fz) + {F^{\textnormal{sub}}(\fz) \over N} \Big) ),
    \end{align*}
    where $\Hat{Z}_{N, \beta}>0$ is a normalizing constant, and \begin{align}
        \label{freedef}  &F_{\beta}(\fz) = -\beta \Tr(\fz^{\intercal} \fJ \fz) + \sum_{i=1}^{m} \sum_{j=1}^{q} \fz^{ij} \log(\fz^{ij}),\\ 
        \nonumber&F^{\textnormal{sub}}(\fz) = {1 \over 2} \sum_{i=1}^{m} \sum_{j=1}^{q} \log(\fz^{ij}) + O(N^{-(q-1)}).
    \end{align}
\end{lemma}
\begin{proof}
    Fix $\fz \in \cPm.$ Since there are $\prod_{i=1}^{m} \frac{N!}{(N \fz^{i1})! \cdots (N\fz^{iq})!}$ configurations with proportion $\fz,$ we get \begin{align*}
        \nu_{N, \beta}(\fz) &= \frac{1}{Z_{N, \beta}} \sum_{\sigma \in \Sigma: \fS(\sigma)= \fz} \exp( \beta N \Tr(\fz^{\intercal} \fJ \fz)) \\
        &= {1 \over Z_{N, \beta}} \prod_{i=1}^{m} \frac{N!}{(N \fz^{i1})! \cdots (N\fz^{iq})!} \exp( \beta N \Tr(\fz^{\intercal} \fJ \fz)).
    \end{align*}
    By Stirling's approximation, \begin{align*}
        \frac{N!}{(N \fz^{i1})! \cdots (N\fz^{iq})!} \approx \frac{1}{(2\pi N)^{(q-1)/2}} \exp( -N \sum_{j=1}^{q} \fz^{ij} \log( \fz^{ij}) - {1 \over 2} \sum_{j=1}^{q} \log \fz^{ij}).
    \end{align*} Therefore, we have \begin{align*}
        \nu_{N, \beta}(\fz) = {1 \over \hat{Z}_{N, \beta}} \exp( -N \Big(F_{\beta}(\fz) + {F^{\operatorname{sub}}(\fz) \over N}\Big)),
    \end{align*}
    where $\hat{Z}_{N, \beta}$ is some normalizing constant, and $F_{\beta}, \ F^{\operatorname{sub}}$ are defined as above.
\end{proof}

 The dominant term $F_{\beta}$ is called \textit{free energy}, which the Glauber dynamics tend to minimize over time. (See \cite{ellisentropy, path, LDPblock}.) We now compute the conditional probability of updating a configuration $\sigma$ at a chosen vertex $v \in V_{i}$ to a specific spin $j \in [q].$ 
\begin{lemma} \label{update lem}
   The probability of updating configuration $\sigma$ by changing to a spin $j \in [q]$ at a chosen vertex $v \in V_{i}$ is given by
    \begin{align} \label{update prob}
    q^{i}_{v, \sigma} (e^{j}) = g_{\beta}^{i j}\big(\fS(\sigma)\big) + O(N^{-1}),
    \end{align}
    where for $i \in [m], \ j \in [q],$ and $\mathbf{z} \in \mathcal{S},$ \begin{align} \label{g def}
    g_{\beta}^{ij}(\mathbf{z})={\exp({2\beta}(\mathbf{J}\mathbf{z})^{ij}) \over \sum_{k=1}^{q} \exp({2\beta }(\mathbf{J}\mathbf{z})^{ik})}.
\end{align}
\end{lemma}
\begin{proof}
    Suppose $\sigma \in \Sigma, \ i \in [m], \ v \in [N], \ \text{and } j \in [q]$ are given. Denote $\sigma^{i}_{v, e^{j}} \in \Sigma,$ obtained by the update on configuration $\sigma$ to have a spin $j$ at $v \in V_{i}.$ By \eqref{Gibbs expression},  \begin{align*}
        q^{i}_{v, \sigma}(e^{j}) = \frac{\exp \big({\beta N } \operatorname{Tr}((\fS(\sigma^{i}_{v, e^j}))^{\intercal}\fJ \fS(\sigma^{i}_{v, e^j})) \big)}{\sum_{l=1}^{q} \exp \big({\beta N } \operatorname{Tr}((\fS(\sigma^{i}_{v, e^l}))^{\intercal}\fJ \fS(\sigma^{i}_{v, e^l})) \big)}.
    \end{align*}
   Let the spin at $v$ be $j_0.$ Then, \begin{align*}
        \sigma^{i}_{v, e^{j}}(w)= 
        \begin{cases}
            j & \text{if } w=v,\\
            \sigma(w) & \text{if } w \neq v.
        \end{cases}
    \end{align*}

    Assume $j \neq j_{0}.$
    We have for each $l \in [m], \ k \in [q],$ \begin{align} \label{sigma' prop}
        \fS(\sigma^{i}_{v, e^{j}})^{lk} = \begin{cases}
            \fS(\sigma)^{ij_0} - \frac{1}{N} & \text{if } l=i, k = j_0, \\
            \fS(\sigma)^{ij} + \frac{1}{N} & \text{if } l=i, k = j, \\
            \fS(\sigma)^{lk} & \text{otherwise.}
        \end{cases}
    \end{align} 
    Therefore, we compute \begin{align*}
        &\Tr(\fS(\sigma^{i}_{v, e^{j}})^{\intercal} \fJ \fS(\sigma^{i}_{v, e^{j}})) = \sum_{k \in [q]} \sum_{l, s \in [m]} \fS(\sigma^{i}_{v, e^{j}})^{lk} \fJ^{ls} \fS(\sigma^{i}_{v, e^{j}})^{sk}\\
        &= \sum_{k \neq j_0, j} \sum_{l, s \in [m]} \fS(\sigma^{i}_{v, e^{j}})^{lk} \fJ^{ls} \fS(\sigma^{i}_{v, e^{j}})^{sk}+ \underbrace{\sum_{l, s \in [m]} \fS(\sigma^{i}_{v, e^{j}})^{lj_0} \fJ^{ls} \fS(\sigma^{i}_{v, e^{j}})^{sj_0}}_{\text(a)}+\underbrace{\sum_{l, s \in [m]} \fS(\sigma^{i}_{v, e^{j}})^{lj} \fJ^{ls} \fS(\sigma^{i}_{v, e^{j}})^{sj}}_{\text(b)}.
    \end{align*}
    Firstly, we expand $\text(a)$ into $\fS(\sigma)$ terms as
    \begin{align*}
        \text(a) &= \sum_{l, s \neq i} \fS(\sigma)^{lj_0} \fJ^{ls} \fS(\sigma)^{sj_0}+2 \sum_{l \neq i} \Big(\fS(\sigma)^{ij_0} - {1 \over N}\Big) \fJ^{li} \fS(\sigma)^{lj_0} + \fJ^{ii} \Big(\fS(\sigma)^{ij_0} - {1 \over N}\Big)^2\\
        &=\sum_{l, s \in [m]} \fS(\sigma)^{lj_0} \fJ^{ls} \fS(\sigma)^{sj_0} - \frac{2}{N} (\fJ \fS(\sigma))^{ij_0} + { \fJ^{ii} \over N^{2}}.
    \end{align*}
    Similarly, we expand the rest of terms as \begin{align*}
        &\text(b)=\sum_{l, s \in [m]} \fS(\sigma)^{lj} \fJ^{ls} \fS(\sigma)^{sj} + \frac{2}{N} (\fJ \fS(\sigma))^{ij} + {\fJ^{ii} \over N^{2}},\\
        &\sum_{k \neq j_0, j} \sum_{l, s \in [m]} \fS(\sigma^{i}_{v, e^{j}})^{lk} \fJ^{ls} \fS(\sigma^{i}_{v, e^{j}})^{sk} = \sum_{k \neq j_0, j} \sum_{l, s \in [m]} \fS(\sigma)^{lk} \fJ^{ls} \fS(\sigma)^{sk},
    \end{align*} where the second line is from ~(\ref{sigma' prop}). This leads us to \begin{align*}
        \Tr(\fS(\sigma^{i}_{v, e^{j}})^{\intercal} \fJ \fS(\sigma^{i}_{v, e^{j}})) = \Tr(\fS(\sigma)^{\intercal} \fJ \fS(\sigma)) + \frac{2}{N} \big((\fJ \fS(\sigma))^{ij} - (\fJ \fS(\sigma))^{ij_0}\big) + \frac{2 \fJ^{ii}}{N^{2}}.
    \end{align*}
Therefore, we obtain \begin{align*}
    &\exp \big({\beta N } \operatorname{Tr}((\fS(\sigma^{i}_{v, e^j}))^{\intercal}\fJ \fS(\sigma^{i}_{v, e^j})) \big) \\
    &=\begin{cases}
        \exp(\beta N \Tr(\fS(\sigma)^{\intercal} \fJ \fS(\sigma))) & \text{if } j = j_0,\\
        \exp \Big(\beta N \Tr(\fS(\sigma)^{\intercal} \fJ \fS(\sigma)) + 2\beta \big( (\fJ \fS(\sigma))^{ij}- (\fJ \fS(\sigma))^{ij_0} \big) + {2\beta \over N} \fJ^{ii} \Big) & \text{if } j \neq j_0,
    \end{cases}
\end{align*}
Using $\exp({{2\beta \over N} \fJ^{ii}}) = 1 + {2 \beta \over N} \fJ^{ii} + O(N^{-2})$, we conclude that \begin{align} \label{q accurate}
q^{i}_{v, \sigma}(e^{j})= \begin{cases}
    g^{ij_0}_{\beta}(\fS(\sigma)) - {2 \beta \over N} \fJ^{ii} g^{ij_0}_{\beta}(\fS(\sigma)) \big(1 - g^{ij_0}_{\beta}(\fS(\sigma)) \big) + O(N^{-2}) & \text{if } j = j_0,\\
    g^{ij}_{\beta}(\fS(\sigma)) + {2 \beta \over N} \fJ^{ii} g^{ij}_{\beta}(\fS(\sigma)) g^{ij_0}_{\beta}(\fS(\sigma)) + O(N^{-2}) & \text{if } j \neq j_0.
\end{cases}
\end{align}
\end{proof}

\begin{remark}
    Up to order $O(1)$, the update probability \eqref{update prob} in this lemma is independent of the choice of vertex $v$ among $V_i = \{1, \cdots, N\}.$ 
\end{remark}

We conclude this section with basic properties of the interaction matrix \(\fJ=J\,\bm\Lambda\).  The proof is by direct computation.

\begin{lemma} \label{matrix property}
Let \(\lambda\in(0,1)\).  Then the following hold for $m \times m$ matrix \begin{align*}
    \bm{\Lambda} = \frac{1}{1+(m-1)\lambda} \begin{bmatrix}
       1      & \lambda & \cdots & \lambda\\
       \lambda & 1      & \cdots & \lambda\\
       \vdots & \vdots & \ddots & \vdots\\
       \lambda & \lambda & \cdots & 1
     \end{bmatrix}
\end{align*}
\begin{enumerate}
  \item $\bm\Lambda$ is positive definite with eigenvalues
  \[
    \frac{1-\lambda}{1+(m-1)\lambda},
    \quad
    1,
  \]
  of multiplicity $m-1$ and $1,$ respectively.
  \item Its operator norms satisfy
  \[
    \lVert\bm\Lambda\rVert_{1\to1}
    = \lVert\bm\Lambda\rVert_{2\to2}
    = 1.
  \]
\end{enumerate}
\end{lemma}

 \medskip

 \subsection{Free energy in the Curie-Weiss-Potts model} We recall the single component case $m=1$ with $J=1$ and $q \geq 3.$ On the simplex $$\cP = \left\{ z=(z_1, \cdots, z_q): z_{i} \geq 0, \ \sum_{i=1}^{q} z_{i} = 1\right\},$$ the free energy is \begin{align*}
    F^{\operatorname{CWP}}_{\beta}(z) = - \beta \| z \|_{2}^{2} + \sum_{j=1}^{q} z^j \log(z^j),
 \end{align*} by a similar computation in \eqref{freedef}. \cite{ellislimit} showed that there is a threshold $\beta_s$ as, \begin{align*}
     \beta_{s}(q)&= \sup \left\{ \beta>0: F^{\CWP}_{\beta} \text{ has a unique critical point } e=(1/q, \cdots, 1/q) \text{ in } \cP \right\}\\
 &= \sup{\left\{\beta >0: \big(1+(q-1)e^{2\beta{1-qx \over q-1}}\big)^{-1}-x \neq 0 \ \textnormal{for all} \ x\in (1/q, 1) \right\}}.
\end{align*}
Furthermore, \cite{Glauber} analyzed the one-variable function \begin{align*}
    d_{\beta}(x) = -x + \big(1+(q-1)e^{2\beta{1-qx \over q-1}}\big)^{-1},
\end{align*}
and proved the following.

\begin{lemma}\cite[Proposition 3.1.]{Glauber} \label{solution study}
    The number of solutions of $d_{\beta}(x) =0$ in $[1/q, 1]$ is as follows. \begin{enumerate}
        \item For $\beta< \beta_{s}(q),$ $d_{\beta}$ has the unique solution $x=1/q.$ In fact, \[
        d_{\beta}(x) < 0 , \quad \text{for all } x \in (1/q, 1).
        \]
        \item For $\beta= \beta_{s}(q),$ $d_{\beta}$ has two solutions $x=1/q$ and $x = x^{*}$ in $(0, 1)$ such that $${d \over dx} d_{\beta}(x^{*}) = 0.$$
        \item For $\beta> \beta_{s}(q),$ $d_{\beta}$ has three solutions $1/q < x_{1}(\beta) < x_{2}(\beta)$ in $[1/q, 1).$ Their derivatives satisfy \[
        {d \over dx} d_{\beta}(1/q) < 0, \ {d \over dx} d_{\beta}(x_{1}) > 0 , \ \text{and } \ {d \over dx} d_{\beta}(x_{2}) < 0.
        \]
    \end{enumerate}
\end{lemma}

On the other hand, \cite{rapid} characterized $\beta_{s}(q)$ via the map $\cG_{\beta}: \cP \to \cP$ by 
\begin{align} \label{rapid beta_s}
    \beta_s(q) &= \sup{\left\{\beta >0: \mathcal{G}^{j}_{\beta}(z) < z^{j} \ \textnormal{for all} \ z \in \cP \ \text{with } z^{j} > 1/q \right\}},
\end{align}
where $\mathcal{G}_{\beta}= (\cG_{\beta}^{1}, \cdots, \cG_{\beta}^{q})$ is \[
 \cG_{\beta}^{j}(z) = \frac{\exp(2 \beta z^{j})}{\sum_{k=1}^{q} \exp(2\beta z^{k})}, \quad \forall j \in [q].
\]
One can check \eqref{rapid beta_s} is equivalent to \eqref{beta_s def} via Jensen's inequality, since \begin{align}
    \nonumber-z^{j} + \cG_{\beta}(z) &= -z^{j} + \frac{1}{1 + \sum_{k \neq j}^{q} \exp(2\beta(z^{k} -z^{j})) }\\
    \nonumber&\leq - z^{j}  + \frac{1}{1 +(q-1) \exp({2\beta \over q-1} \sum_{k \neq j} ^{q} (z^{k} -z^{j})) }\\
    \label{g and d} &= - z^{j}  + \frac{1}{1 +(q-1) \exp({2\beta \over q-1}(1 -q z^{j})) } = d_{\beta}(z^{j}).
\end{align}

\begin{remark}
The function $g_{\beta}^{ij}$ in \eqref{g def} can be written as \begin{align*}
        g^{ij}_{\beta}(\fz) = \cG_{\beta}^{j}((\fJ \fz)^{i*}), \quad \fz \in \cPm.
    \end{align*}
\end{remark}

\medskip

\subsection{Large deviation results}
To analyze mixing behavior, we examine the structure of the free energy $F_{\beta}$ on the proportion space $\cPm.$ When $N$ is large, the empirical proportion $\fS(\sigma)$ concentrates near global minimizers of $F_{\beta}.$ \cite{LDPblock} showed that $\fS(\sigma)$ satisfies a \textit{large deviation principle}.

\begin{defn}
    A sequence of probability measures $(\nu_{N, \beta})_{N \in \bN}$ on $\cPm$ satisfies a large deviation principle with rate function $I_{\beta}: \cPm \to [0, \infty],$ if for every closed set $C$ and open set $U$ in $\cPm,$ \begin{align}
     \label{limsup}  &\limsup_{N \to \infty}{1 \over N}\log\big( \nu_{N, \beta}[\fz \in C]\big) \leq - \inf_{\fz \in C }I_{\beta}(\mathbf{z}),\\
    \nonumber&\liminf_{N \to \infty}{1 \over N}\log\big( \nu_{N, \beta} [\fz \in U ]\big) \geq - \inf_{\fz \in U }I_{\beta}(\mathbf{z}).
\end{align}
\end{defn}

\begin{prop}\cite[Theorem 3.2.]{LDPblock} \label{LDP thrm}The measure $\nu_{N, \beta}$ satisfies a large deviation principle with rate function $I_{\beta}$, \begin{align*}
    I_{\beta}(\fz) = F_{\beta}(\fz) - \min_{\fw \in \cPm} F_{\beta}(\fw),
\end{align*}
where $F_{\beta}$ is defined in \eqref{freedef}.
\end{prop}

Define the set of global minimizers \begin{align*}
    \mathcal{E}_{\beta}= \argmin_{z \in \cPm} F_{\beta}(\fz) = \{ \fz \in \cPm: I_{\beta}(\fz) = 0\}. 
\end{align*}

\begin{prop} \cite[Theorem 4.7.]{LDPblock} \label{phase trans}In the homogeneous multi-component Potts model, there is a phase transition at $\beta = {\beta_{c}(q) / J}.$ Specifically, \begin{align*}
        \mathcal{E}_{\beta} = \begin{cases}
            \{\fe\} & \text{if } \beta< {\beta_{c}(q) / J}, \\
            \{\fe, \ \bm{\eta}^{j}(\beta): \ 1 \leq j \leq q\} & \text{if } \beta= {\beta_{c}(q) / J},\\
            \{\bm{\eta}^{j}(\beta): \ 1 \leq j \leq q\} & \text{if } \beta > {\beta_{c}(q) / J},
        \end{cases}
    \end{align*}
    where \begin{align*}
        \beta_{c}(q) = {(q-1) \over q-2} \log(q-1). 
    \end{align*}
    $\fe$ is $m \times q$ matrix with all entries $q^{-1},$ and  
    $\bm{\eta}^{1}(\beta)$ is $m \times 
    q$ matrix as
    \begin{align*}
        \bm{\eta}^{1}(\beta) =
       \begin{bmatrix}
1-(q-1) t(\beta J) \ \; \ t(\beta J) \; \ t(\beta J) \; \ \cdots \ \; \ t(\beta J)\\
1-(q-1) t(\beta J) \ \; \ t(\beta J) \; \ t(\beta J) \; \ \cdots \ \; \ t(\beta J)\\
\vdots\\
1-(q-1) t(\beta J) \ \; \ t(\beta J) \; \ t(\beta J) \; \ \cdots \ \; \ t(\beta J)
\end{bmatrix}.
    \end{align*}
    For $2 \leq j \leq q$, the $\bm{\eta}^{j}(\beta)$ is obtained by swapping the first column of $\bm{\eta}^{1}(\beta)$ with its $j$-th column. The $t = t(\beta J)$ is the smallest positive solution of the equation \begin{align*}
        t = \frac{1}{q-1 + \exp(2\beta J (1-qt))}.
    \end{align*} 
\end{prop}

Proposition ~\ref{phase trans} characterizes the discontinuous phase transition in our model, by describing how the set of global minimizers of $F_{\beta}$ changes. When $\beta<\beta_{c}(q)/J$, the system is in the disordered phase with the unique global minimizer $\fe.$ As $\beta \geq \beta_{c}(q)/J$, new global minimizers $\bm{\eta}^{j}(\beta)$ emerge.

\bigskip
\section{Fast mixing in high-temperature $\beta< \beta_{s}(q)/J$} \label{high} 
In this section, we establish an $O(N \log N)$ bound when $\beta< \beta_{s}(q)/J$. Path coupling is a standard technique in probability theory that quantifies the  contraction in coupling distance between two nearby configurations via a path metric. However, it requires the initial configurations to be neighbors. To overcome this limitation, \cite{blumepath} introduced the \textit{aggregate path} approach, which connects two configurations by a sequence of neighboring configurations. This method preserves the strength of path coupling argument while allowing it to apply to configurations that are not necessarily adjacent.  

\subsection{Greedy Coupling} We select a vertex $v$ uniformly at random and update both configurations at $v$ according to their conditional Gibbs measures. We first define the \emph{Hamming distance} between two configurations $\sigma,\tau\in\Sigma$.

\begin{defn}
    For configurations $\sigma=(\sigma^1, \sigma^2, \cdots, \sigma^m)$ and $\tau=(\tau^1, \tau^2, \cdots, \tau^m)$, the Hamming distance $d$ is defined by \begin{align*}
        d(\sigma, \tau)=\sum_{i=1}^{m} \sum_{v \in V_i} \mathbbm{1}_{\{\sigma(v) \neq \tau(v)\}}.
    \end{align*}
\end{defn}

We construct a \textbf{greedy coupling} of two copies of Glauber dynamics. At each time step $t$, let $\sigma_t=(\sigma^1_t, \cdots, \sigma^m_t)$ and $\tau_t = (\tau^1_t, \cdots, \tau^m_t).$  \begin{itemize}
    \item[(1)] Choose a vertex $v$ uniformly at random. Suppose $v \in V_{i}$ for some $i \in [m].$
    \item[(2)] Draw spins $J_{t+1}, \tilde{J}_{t+1},$ via the maximal coupling of $q^{i}_{v, \sigma_{t}}(\cdot)$ and  $q^{i}_{v, \tau_{t}}(\cdot).$
    \item[(3)] Update the spin at $v$ to $J_{t+1}$ in $\sigma_t$ and $\tilde{J}_{t+1}$ in $\tau_t$.
    \item[(4)] Remain the rest of spins to be unchanged at the other vertices $w \neq v.$
\end{itemize}
The maximal coupling synchronizes the two processes as much as possible. (See \cite{markov}.)

\begin{lemma}\cite [Proposition 4.7]{markov}
    Let $\mu$ and $\nu$ be two probability measures on a space $\Sigma.$ Then, \begin{align*}
    \| \mu - \nu \|_{\operatorname{TV}} = \inf \{ \mathbb{P}[X \neq Y]: (X, Y) \text{ is a coupling of } \mu \text{ and } \nu \},
    \end{align*}
    and there exists a coupling achieving the infimum called, the {maximal coupling}.
\end{lemma}

We denote the measure and expectation under this greedy coupling by $\mathbb{P}^{GC}$ and $\mathbb{E}^{GC}.$ We next compute the expected drift of $d(\sigma_t, \tau_t).$ 

\begin{lemma} \label{greedy coupling}
    We consider one step of the greedy coupling, updating $\sigma=(\sigma^1, \sigma^2, \cdots, \sigma^m)$ and $\tau=(\tau^1, \tau^2, \cdots, \tau^m)$ to $\sigma'$ and $\tau'.$ Then, we have \begin{align*}
        \mathbb{E}^{GC} \big[d(\sigma', \tau') \big] = \big[1-{1 \over mN} (1-\Gamma)\big] d(\sigma, \tau),
    \end{align*} where the contraction constant $\Gamma= \Gamma(\sigma, \tau)$ satisfies \begin{align} \label{Gamma def}
        \Gamma \leq  \frac{ \sum_{i=1}^{m} \sum_{j=1}^{q}| g^{ij}_{\beta}(\fS(\sigma)) - g^{ij}_{\beta}(\fS(\tau))| + O(N^{-1})} {\| \fS(\sigma) - \fS(\tau) \|_{1}}
    \end{align}
\end{lemma}
\begin{proof}
By the construction of the greedy coupling and Lemma \ref{update lem}, for each $i \in [m]$, \begin{align*}
    \mathbb{P}^{GC}[ J_{1} \neq \tilde{J}_{1} \big| v \in V_i] &= \| q^{i}_{v, \sigma} - q^{i}_{v, \tau} \|_{\TV}= {1 \over 2} \sum_{j=1}^{q} | g^{ij}_{\beta}(\fS(\sigma)) - g^{ij}_{\beta}(\fS(\tau))| + O(N^{-1}).
\end{align*} Since this probability is independent of the choice of $v \in V_{i},$ set $\kappa^{i} =   \| q^{i}_{v, \sigma} - q^{i}_{v, \tau} \|_{\TV}$ and \begin{align*}
    d^{i} = \sum_{v \in V_{i}}  \mathbbm{1}_{\{\sigma(v) \neq \tau(v)\}},
\end{align*} so that $d(\sigma, \tau) = \sum_{i=1}^{m} d^{i}.$ Then, we have \begin{align*}
    \mathbb{E}^{GC} \big[d(\sigma', \tau') \big] &= d(\sigma, \tau) + {1 \over m} \sum_{i=1}^{m} \Big[{N-d^{i} \over N} \kappa^{i} - {d^{i} \over N} (1-\kappa^{i}) \Big]\\
    \nonumber&= \Big(1 - {1 \over mN} \Big) d(\sigma, \tau) + {1 \over m} \sum_{i=1}^{m} \kappa^{i} \\
    &= d(\sigma, \tau) \Bigg[1 - {1 \over mN} \Big( 1 - \frac{ \sum_{i=1}^{m} \kappa^{i} }{2d(\sigma, \tau)/N} \Big) \Bigg]\\
    &\leq d(\sigma, \tau) \Bigg[1 - {1 \over mN} \Big( 1 - \frac{ \sum_{i=1}^{m} | g^{ij}_{\beta}(\fS(\sigma)) - g^{ij}_{\beta}(\fS(\tau))| + O(N^{-1}) }{\| \fS(\sigma) - \fS(\tau)\|_{1}} \Big) \Bigg],
\end{align*}
where the inequality is from \begin{align*}
    \frac{2 d(\sigma^{i}, \tau^{i})}{N} \geq \| \fS^{i*}(\sigma) - \fS^{i*}(\tau) \|_{1}, \quad \text{for each } i \in [m].
\end{align*}
\end{proof}
\medskip


\subsection{Aggregate path method}Our goal is to show that the expected mean coupling distance is contracted by aggregating intermediate distances along a monotone path in the configuration space. We first define a monotone path.

\begin{defn}
    Let $\sigma, \tau \in \Sigma$ be two configurations. A path \begin{align*}
        \pi: \sigma=\omega_0, \ \omega_1, \ \cdots, \ \omega_r = \tau
    \end{align*}in $\Sigma$ is called a \textit{monotone} if it satisfies \begin{enumerate}
        \item $d(\sigma, \tau)=\sum_{s=1}^{r}d(\omega_{s-1}, \omega_{s}),$
        \item For each $(i, j) \in [m] \times [q]$ , $\fS(\omega_{s})^{ij}$ is monotonic in $s=0, 1, \cdots, r.$
    \end{enumerate}
\end{defn} 

To construct a monotone path in $\Sigma$, we first construct its analogue in the proportion space $\cPm.$
\begin{remark} \label{path correspondence}
    As in \cite[page 76]{path}, one checks that if there is a path \begin{align*}
        \pi^{\operatorname{prop}} : \fS(\sigma) = \fz_{0}, \ \fz_{1}, \ \cdots, \ \fz_{r} = \fS(\tau), \quad \sigma, \tau \in \Sigma,
    \end{align*}
    in $\cPm,$ monotone in each coordinate $(i, j) \in [m] \times [q],$ then there exists the corresponding monotone path \begin{align*}
        \pi : \sigma = \omega_{0}, \ \omega_{1}, \ \cdots, \ \omega_{r} = \tau,
    \end{align*} in $\Sigma$
    such that $\fS(\omega_s) = \fz_{s}$ for all $0 \leq s \leq r.$
\end{remark}

Assuming such a monotone path exists, we can decompose the contraction constant \eqref{Gamma def} into terms along each step of the path.

\begin{lemma} \label{contraction under monotone}
    Under the same hypotheses of Lemma \ref{greedy coupling}, assume that for sufficiently small $\eps>0,$ there is a coordinate-wise monotone path \begin{align*}
        \pi^{\operatorname{prop}}: \fS(\sigma) = \fz_{0}, \ \fz_{1}, \ \cdots, \ \fz_{r} = \fS(\tau)
    \end{align*} in $\cPm$
    with $\| \fz_{s-1} - \fz_{s} \|_{1} \in (\eps, 2\eps], \ \text{for each } s \in [1, r].$ Then, the contraction constant $\Gamma$ satisfies \begin{align*}
        \Gamma \leq \frac{\sum_{s=1}^{r} \sum_{i=1}^{m} \sum_{j=1}^{q} | \inner{ \nabla g^{ij}_{\beta} (\fz_{s-1})}{\fz_{s-1} - \fz_{s}}| + O(\eps^{2}) + O(\eps N^{-1}) + O(N^{-2})}{\| \fS(\sigma) - \fS(\tau) \|_{1}}.
    \end{align*} 
\end{lemma}

\begin{proof}
    By Lemma \ref{greedy coupling}, \begin{align*}
        \Gamma = \frac{ \sum_{i=1}^{m} \kappa^{i}}{d(\sigma, \tau) / N} = \frac{1}{2d(\sigma, \tau)/N} \sum_{i=1}^{m} \sum_{j=1}^{q} |q^{i}_{v, \sigma}(e^{j}) - q^{i}_{v, \tau}(e^{j})|.
    \end{align*}

    First suppose $\| \fS(\sigma) - \fS(\tau) \|_{1} \in (\eps, 2 \eps ]$ for some small $\eps>0.$

    \medskip
    If $\sigma(v) = \tau(v) = j_{0}$ with some $j_0 \in [q],$ then by \eqref{q accurate}, one shows \begin{align*}
        q^{i}_{v, \sigma}(e^j) - q^{i}_{v, \tau}(e^j) = & g^{ij}_{\beta}(\fS(\sigma)) - g^{ij}_{\beta}(\fS(\tau)) + O(N^{-2})\\ + &{2 \beta \fJ^{ii} \over N}  \begin{cases}
            \nonumber g^{ij}_{\beta}(\fS(\sigma)) g^{ij_0}_{\beta}(\fS(\sigma)) - g^{ij}_{\beta}(\fS(\tau))g^{ij_0}_{\beta}(\fS(\tau)) & \text{if } j \neq j_0\\
            -\nonumber\big[g^{ij}_{\beta}(\fS(\sigma)) (1 -g^{ij}_{\beta}(\fS(\sigma))) - g^{ij}_{\beta}(\fS(\tau))(1 - g^{ij}_{\beta}(\fS(\tau)))\big] & \text{if } j = j_0\\
        \end{cases}\\
        = & g^{ij}_{\beta}(\fS(\sigma)) - g^{ij}_{\beta}(\fS(\tau)) + O\Big({\eps \over N}\Big) + O(N^{-2}),
    \end{align*}
    since $g^{ij}_{\beta}g^{ij_0}_{\beta}$ and $g^{ij}_{\beta} (1- g^{ij}_{\beta})$ are Lipschitz continuous for all $i \in [m], \ j, j_0 \in [q].$ Applying $\| \fS(\sigma) - \fS(\tau)\|_{1} = O(\eps),$ we deduce \begin{align}
        \nonumber &\sum_{i=1}^{m} \sum_{j=1}^{q} |q^{i}_{v, \sigma}(e^{j}) - q^{i}_{v, \tau}(e^{j})|\\ \nonumber&= \sum_{i=1}^{m} \sum_{j=1}^{q} |g^{ij}_{\beta}(\fS(\sigma)) - g^{ij}_{\beta}(\fS(\tau))| + O\Big( {\eps \over N} \Big) + O(N^{-2})\\
        \label{q approx} &=\sum_{i=1}^{m} \sum_{j=1}^{q} |\inner{\nabla g^{ij}_{\beta}(\fS(\sigma))}{\fS(\sigma) - \fS(\tau)}| + O(\eps^{2}) + O\Big( {\eps \over N} \Big) + O(N^{-2}).
    \end{align}

    In a more general situation when $\sigma(v) \neq \tau(v),$ the same estimate holds as in \eqref{q approx}. Now let \begin{align*}
        \pi^{\operatorname{prop}}: \fS(\sigma) = \fz_{0}, \ \fz_{1}, \ \cdots, \ \fz_{r} = \fS(\tau)
    \end{align*} be the given monotone path in $\cPm$ such that $\| \fz_{s-1} - \fz_{s}\|_{1} \in (\eps, 2\eps].$
    By Remark \ref{path correspondence}, it lifts to a corresponding monotone path \begin{align*}
        \pi: \sigma = \omega_{0}, \ \omega_{1}, \ \cdots, \ \omega_{r} = \tau  \ \ \text{on } \Sigma \ \ \text{with } \fS(\omega_s) = \fz_{s}.
    \end{align*}
    Summing \eqref{q approx} along each step, we obtain \begin{align*}
        \sum_{i=1}^{m} \sum_{j=1}^{q} |q^{i}_{v, \sigma}(e^{j}) - q^{i}_{v, \tau}(e^{j})| &\leq \sum_{i=1}^{m} \sum_{j=1}^{q} \sum_{s=1}^{r} |q^{i}_{v, \omega_{s-1}}(e^{j}) - q^{i}_{v, \omega_{s}}(e^{j})| \\
        &= \sum_{i=1}^{m} \sum_{j=1}^{q} \sum_{s=1}^{r} | \inner{\nabla g^{ij}_{\beta}( \fz_{s-1})}{\fz_{s-1} - \fz_{s}}| + O(\eps^{2}) + O\Big( {\eps \over N} \Big) + O(N^{-2}),
    \end{align*}
    which is what we wanted.
\end{proof}

To prove contraction in the high-temperature regime $\beta< \beta_s(q) / J$, we want \begin{align*}
    \frac{\sum_{i=1}^{m} \sum_{j=1}^{q} |g_{\beta}^{ij} (\fS(\sigma)) - g_{\beta}^{ij}(\fS(\tau))|}{\| \fS(\sigma) - \fS(\tau)\|_{1}} < 1, \quad \text{uniformly in } \sigma,
\end{align*} 
and $\tau$ with the necessary conditions. We first investigate this ratio when both $\sigma$ and $\tau$ lie close to the equilibrium $\fe=q^{-1} \mathbf{1}_{m \times q}$.

\begin{lemma}
    We have \begin{align*}
        \lim_{{\fz \to \fe, \ \fz \in \cPm}} \frac{\sum_{i=1}^{m} \sum_{j=1}^{q} |g^{ij}_{\beta}(\fz) - g^{ij}_{\beta}(\fe)|}{\| \fz - \fe\|_{1}} \leq {2\beta J \over q}.
    \end{align*}
    Consequently, if $\beta < \beta_{s}(q)/J,$ then this limit is strictly less than $1.$
\end{lemma}
\begin{proof}
    Recall that \begin{align*}
        g^{ij}_{\beta}(\fz) = \frac{ \exp(2\beta (\fJ \fz)^{ij})}{\sum_{k=1}^{q} \exp(2\beta (\fJ \fz)^{ik})}.
    \end{align*} We compare $\exp(2\beta (\fJ \fz)^{ij})$ to $\exp(2\beta (\fJ \fe)^{ij})$ via the Taylor expansion $\exp(t) = 1+ t + O(t^{2}).$ Writing $\bdel = \fz - \fe,$ we have \begin{align*}
        \exp(2\beta (\fJ \fz)^{ij} - 2\beta (\fJ \fe)^{ij}) &= \prod_{u=1}^{m} \exp( 2\beta \fJ^{iu} \bdel^{uj} )\\
        &= \prod_{u=1}^{m} \Big[1 + 2\beta \fJ^{iu} \bdel^{uj} + O(| \bdel^{uj}|^{2}) \Big]\\
        &= 1 + 2\beta (\fJ\bdel)^{ij} + O(\| \bdel^{*j}\|_{2}^{2}).
    \end{align*}
    Since $\fJ \fe = \fe,$ it follows that \begin{align*}
         \exp(2\beta (\fJ \fz)^{ij}) &= \exp(2\beta / q) [ 1 + 2\beta (\fJ\bdel)^{ij} + O(\| \bdel^{*j}\|_{2}^{2})], \\
         \sum_{k=1}^{q}  \exp(2\beta (\fJ \fz)^{ik}) &= \exp(2\beta / q) [q + O(\| \bdel \|_{2}^{2})].
    \end{align*}
    Combining these two, we obtain \begin{align*}
        g^{ij}_{\beta}(\fz) - g^{ij}_{\beta}(\fe) = g^{ij}_{\beta}(\fz) - {1 \over q} ={2 \beta \over q} (\fJ \bdel)^{ij} (1 + O(\|\bdel\|_{2})).
    \end{align*}
    Summing over $i, \ j$ gives \begin{align*}
        \sum_{i=1}^{m} \sum_{j=1}^{q} |g^{ij}_{\beta}(\fz) - g^{ij}_{\beta}(\fe)| = {2 \beta \over q} \sum_{i=1}^{m} \sum_{j=1}^{q} |(\fJ \bdel)^{ij}| + O(\|\bdel\|_{2}^{2}).
    \end{align*}
   By Lemma \ref{matrix property},   \begin{align*}
        \sum_{i=1}^{m} \sum_{j=1}^{q} | (\fJ \bdel)^{ij} | = \sum_{j=1}^{q} \| \fJ \bdel^{*j} \|_{1} \leq \| \fJ \|_{1 \to 1} \sum_{j=1}^{q} \| \bdel^{*j}\|_{1} = J \| \bdel\|_{1}.
    \end{align*}
    Thus, \begin{align*}
        \frac{ \sum_{i=1}^{m} \sum_{j=1}^{q} |g^{ij}_{\beta}(\fz) - g^{ij}_{\beta}(\fe)| }{\| \fz - \fe\|_{1}} \leq {2\beta J \over q} + O(\| \fz - \fe\|_{1}),
    \end{align*}
    and the claimed limit follows as $\fz \to \fe.$ 
\end{proof}

Since $\cPm$ is a compact subset of $\bR^{m \times q}$ and each $g^{ij}_{\beta}$ is continuous, we immediately obtain the following. 
\begin{remark} \label{when z, w are close}
    By \cite[Proposition 3.1]{Glauber}, $\beta_{s}(q) < q/2.$ Hence if $\beta < {\beta_{s}(q)/ J},$ there exists a constant $\gamma_0 \in (0, 1)$  such that for any two distinct points $\fz, \ \fw \in \cPm$ with  $\|\fz - \fe\|_{1}, \ \| \fw - \fe \|_{1} < \eps,$ one has \begin{align*}
        \frac{ \sum_{i=1}^{m} \sum_{j=1}^{q} |g^{ij}_{\beta}(\fz) - g^{ij}_{\beta}(\fw)|}{\|\fz - \fw\|_{1}} < 1- \gamma_0,
    \end{align*}
    for sufficiently small $\eps>0.$
\end{remark}
Heuristically, when $\fS(\sigma)$ and $\fS(\tau)$ are close, we approximate \begin{align*}
     \frac{\sum_{i=1}^{m} \sum_{j=1}^{q} |g_{\beta}^{ij} (\fS(\sigma)) - g_{\beta}^{ij}(\fS(\tau))|}{\| \fS(\sigma) - \fS(\tau)\|_{1}} \approx  \frac{\sum_{i=1}^{m} \sum_{j=1}^{q} \int_{\rho} | \inner{\nabla g_{\beta}^{ij}( \fy)}{d\fy}}{\| \fS(\sigma) - \fS(\tau)\|_{1}},
\end{align*} 
where $\rho$ is some continuous monotone path in $\cPm$ linking $\fS(\sigma)$ and $\fS(\tau).$  
\begin{defn} \label{D def}
    For any two points $\fx, \fz \in \cPm$, define the \textit{aggregated $g$-variation}\begin{align*}
        D^{g}(\mathbf{z}, \mathbf{x})= \sum_{i=1}^{m} \sum_{j=1}^{q} \int_{\rho} \Big| \big\langle \nabla g_{\beta}^{ij}(\fy), d\fy \big\rangle \Big|,
    \end{align*} where $\rho$ is the straight line path in $\cPm$ from $\fx$ to $\fz.$
\end{defn}
From now on, unless otherwise specified, we take $\rho$ to be the straight line joining $\fe$ to $\fz.$ 

\begin{lemma} \label{main claim}
    Let $\beta<\beta_{s}(q)/J$ and $\fz \neq \fe$ lie in $\cPm.$ Then, for the straight line path $\rho $ from $\fe$ to $\fz$ in $\cPm,$ we have \begin{align*}
        \frac{D^{g}(\fz, \fe)}{\|\fz - \fe\|_{1}} < 1.
    \end{align*}
\end{lemma}
\begin{proof}
We parameterize the straight line path $\rho$ by \begin{align*}
    \rho^{ij}(t) = {1 \over q}(1-t) + \fz^{ij} t, \quad t \in [0, 1],
\end{align*} for each $(i, j) \in [m] \times [q].$ Then, \begin{align} \label{D cal}
D^{g}(\mathbf{z}, \fe)&= \sum_{i=1}^{m} \sum_{j=1}^{q} \int_{\rho} \Big| \big\langle \nabla g_{\beta}^{ij
}(\fy), d\fy \big\rangle \Big|
= \sum_{i=1}^{m} \sum_{j=1}^{q} \int_{0}^{1} \Big| {d \over dt} g_{\beta}^{ij}(\rho(t)) \Big| dt.
\end{align} Since  \begin{align*}
    {d \over dt} \exp( 2\beta (\fJ \rho(t))^{ij} ) = 2\beta ( \fJ (\fz - \fe))^{ij} \exp( 2\beta (\fJ \rho(t))^{ij} ),
\end{align*}
it follows that
\begin{align*} 
{d \over dt} g_{\beta}^{ij}(\rho(t))={2 \beta}g_{\beta}^{ij}(\rho(t))\big[(\fJ(\mathbf{z}-\mathbf{e}))^{ij}- \inner{ (\mathbf{J}(\mathbf{z}-\mathbf{e}))^{i*}}{g_{\beta}^{i *}(\rho(t))} \big].
\end{align*} Next, we observe \begin{align*}
    {d \over dt}  \inner{ (\mathbf{J}(\mathbf{z}-\mathbf{e}))^{i*}}{g_{\beta}^{i *}(\rho(t))} &= \sum_{k=1}^{q} (\fJ (\fz - \fe))^{ik} {d \over dt} g^{ik}_{\beta}(\rho(t))\\
    &=2 \beta \Big[ \sum_{k=1}^{q} \big((\fJ(\mathbf{z}-\mathbf{e}))^{ij} \big)^{2} g^{ik}_{\beta}(\rho(t)) - \Big(\sum_{k=1}^{q} (\fJ(\mathbf{z}-\mathbf{e}))^{ij} g^{ik}_{\beta}(\rho(t)) \Big)^{2}  \Big]\\
    &= 2\beta \Var_{g^{i*}_{\beta}(\rho(t))} \big[ (\fJ (\fz - \fe))^{i*} \big] \geq 0,
\end{align*}
and at $t=0,$ \begin{align*}
    \inner{ (\mathbf{J}(\mathbf{z}-\mathbf{e}))^{i*}}{g_{\beta}^{i *}(\rho(0))} = \inner{ (\mathbf{J}(\mathbf{z}-\mathbf{e}))^{i*}}{\fe} = {1 \over q} \sum_{k=1}^{q} (\fJ (\fz - \fe))^{ik} = 0. 
\end{align*}
Therefore, $  \inner{ (\mathbf{J}(\mathbf{z}-\mathbf{e}))^{i*}}{g_{\beta}^{i *}(\rho(t))} \geq 0$ for all $t \in [0, 1],$ and we conclude that \begin{enumerate}
    \item If $(\fJ (\fz - \fe))^{ij} \leq 0,$ then $g_{\beta}^{ij}(\rho(t))$ is monotonically decreasing in $t.$
    \item If $(\fJ (\fz - \fe))^{ij} > 0,$ then there is at most one $t_{*}^{ij}$ such that \begin{enumerate}
        \item \label{case a} $g_{\beta}^{ij}(\rho(t))$ is monotonically increasing for $0\leq t\leq t_{*}^{ij},$
        \item $g_{\beta}^{ij}(\rho(t))$ is decreasing for $t> t_{*}^{ij},$
    \end{enumerate}
    If such $t_{*}^{ij}$ does not exist, we set $t_{*}^{ij} = 1.$
\end{enumerate}
   
For each $i \in [m],$ we define the index set \begin{align*}
\mathcal{A}^{i}({\mathbf{z}})=\big\{j \in [q]: (\fJ (\fz - \fe))^{ij} >0\big\}=\big\{j \in [q]: (\mathbf{J}\mathbf{z})^{ij}> J/q\big\}. \end{align*}
If $j \notin \mathcal{A}^{i}(\fz)$, then  \begin{align} \label{D cal_not in}
\end{align}
Otherwise, if $j \in \mathcal{A}^{i}(\mathbf{z})$, then \begin{align}\label{D cal_in}
\int_{0}^{1}\Big|{d \over dt} g_{\beta}^{ij}(\rho(t))\Big| dt &=  \int_{0}^{t_{*}^{ij}} {d \over dt} g_{\beta}^{ij}(\rho(t)) \ dt - \int_{t_{*}^{ij}}^{1} {d \over dt} g_{\beta}^{ij}(\rho(t)) \ dt  \\
&\nonumber = 2g_{\beta}^{ij}(\rho(t_*^{ij}))-g_{\beta}^{ij}(\mathbf{z})-g_{\beta}^{ij}(\fe).
\end{align}
Plugging ~(\ref{D cal_not in}) and ~(\ref{D cal_in}) into ~(\ref{D cal}), we have \begin{align*}
\nonumber D^{g}(\mathbf{z}, \fe) &= \sum_{i=1}^{m} \Big[ \sum_{j \in \mathcal{A}^{i}(\fz)} \big( 2g_{\beta}^{ij}\big(\rho(t_*^{ij})\big)-g_{\beta}^{ij}(\mathbf{z})-g_{\beta}^{ij}(\fe) \big) + \sum_{j \notin \mathcal{A}^{i}(\fz)} \big(g_{\beta}^{ij}(\fe)-g_{\beta}^{ij}(\mathbf{z}) \big) \Big]\\
 &= 2 \sum_{i=1}^{m} \sum_{j \in \mathcal{A}^{i}(\fz)} \big( g_{\beta}^{ij}(\rho(t_*^{ij}))-g_{\beta}^{ij}(\fe) \big).
\end{align*} 
Finally, for each $j \in \mathcal{A}^{i}(\fz),$ $i \in [m]$, and all $t \in [0, 1],$ \begin{align*}
    (\bm{\Lambda} \rho(t))^{ij} = (1-t) (\bm{\Lambda} \fe)^{ij} + t (\bm{\Lambda} \fz)^{ij} \in \Big[ {1 \over q}, (\bm{\Lambda} \fz)^{ij} \Big]. 
\end{align*} Since  $(\bm{\Lambda} \rho(t_{*}^{ij}) )^{i*} \in \cP$ and by \eqref{g and d}, we have \begin{align*}
    g^{ij}_{\beta}( \rho(t_{*}^{ij})) &\leq d_{\beta J} ( (\bm{\Lambda} \rho(t_*^{ij}))^{ij}) + (\bm{\Lambda} \rho(t_*^{ij}))^{ij}<  (\bm{\Lambda} \rho(t_*^{ij}))^{ij} \leq (\bm{\Lambda} \fz)^{ij}.
\end{align*}
Thus, \begin{align*}
    D^{g} (\fz, \fe) <  \sum_{i=1}^{m} \sum_{j=1}^{q} \big| (\bm{\Lambda}( \fz - \fe))^{ij} \big|
    \leq \| \bm{\Lambda}\|_{1 \to 1} \| \fz - \fe \|_{1} = \| \fz - \fe \|_{1}.
\end{align*}
\end{proof}

\begin{remark} \label{similar D and g sum}
    Since $\rho$ is the straight line from $\fe$ to $\fz,$ we have \begin{align*}
         \lim_{{\fz \to \fe, \ \fz \in \cPm}} \frac{D^{g}(\fz, \fe)}{\| \fz - \fe\|_{1}} =  \lim_{{\fz \to \fe, \ \fz \in \cPm}} \frac{\sum_{i=1}^{m} \sum_{j=1}^{q} |g^{ij}_{\beta}(\fz) - g^{ij}_{\beta}(\fe)|}{\| \fz - \fe\|_{1}}.
    \end{align*}
\end{remark}
We obtain the following stronger uniform bound in $\fz.$

\begin{cor} \label{D bound}
    For $\beta < {\beta_s(q) / J},$ there exists a constant ${\gamma} \in (0, 1)$ such that \begin{align*}
        D^{g}(\fz, \fe) \leq (1 - {\gamma}) \|\fz - \fe\|_{1},
    \end{align*}
    uniformly in $\fz \in \cPm.$
\end{cor}
\begin{proof}
   By Remark \ref{when z, w are close} and Remark \ref{similar D and g sum}, there exists $\gamma_1 \in (0, 1)$ such that \begin{align*}
       D^{g}(\fz, \fe) \leq (1- \gamma_1) \| \fz - \fe\|_{1},
   \end{align*}
   whenever $\| \fz - \fe \|_{1} < \eps_0,$ for some sufficiently small $\eps_0>0.$ Let $K = \{ \fz \in \cPm: \| \fz - \fe \|_{1} \geq \eps_0 \}$. The function \begin{align*}
       \fz \mapsto \mathcal{R}(\fz) = \frac{D^{g}(\fz, \fe)}{\|\fz - \fe \|_{1}}
   \end{align*}
   is continuous on the compact set $K,$ so it attains maximizer at some $\fz(\eps_0) \in K,$ i.e.,  \begin{align*}
       \mathcal{R}(\fz(\eps_0)) = \max_{\fz \in K} \mathcal{R}(\fz).
   \end{align*}
   Applying Lemma \ref{main claim} to $\fz(\eps_0) \neq \fe$, there exists $\delta = \delta(\eps_0) \in (0, 1)$ such that \begin{align*}
       \mathcal{R}(\fz(\eps_0)) = 1 - \delta(\eps_0).
   \end{align*}
   Thus, for all $\fz \in K,$ \begin{align*}
       \frac{D^{g}(\fz, \fe)}{\| \fz - \fe \|_{1}} \leq 1 - \delta(\eps_0).
   \end{align*}
   By choosing $\gamma = \min\{ \gamma_1, \delta(\eps_0)\},$ we obtain the uniform bound \begin{align*}
        D^{g}(\fz, \fe) \leq (1- \gamma) \| \fz - \fe \|_{1} \quad \text{for all } \fz \in \cPm.
   \end{align*}
\end{proof}

\begin{remark} \label{monotone approx}
    Since $\rho$ is straight line from $\fe$ to $\fz,$ its curvature is zero. Therefore, for any sufficiently small $\eps>0$ and any $\fz \in \cPm$ with $\| \fz - \fe\|_{1} > \eps,$ one can construct a discrete coordinate-wise monotone path \begin{align*}
        \pi : \fe = \fz_0, \ \fz_1, \ \cdots, \ \fz_r = \fz
    \end{align*} in $\cPm$
    such that $\| \fz_{s} - \fz_{s-1}\|_{1} \in (\eps, 2\eps]$ for each $1 \leq s \leq r,$ and a constant $C>0$ such that \begin{align*}
        \Big|D^{g}(\fz, \fe) - \sum_{i=1}^{m} \sum_{j=1}^{q} \sum_{s=1}^{r} | \inner{\nabla g^{ij}_{\beta}(\fz_{s-1})}{\fz_{s} - \fz_{s-1}}| \Big| \leq Cr \eps^{2}.
    \end{align*}
    (See \cite[Property 6.5]{path}.)
\end{remark}

Now we prove the main contraction.

\begin{prop} \label{main contraction}
    Let $\beta < {\beta_s(q) / J}.$ According to the greedy coupling, we update $(\sigma, \tau)$ to $(\sigma', \tau') \in \Sigma \times \Sigma$ by one time step. Then, there exists uniform constants $\alpha \in (0, 1)$ and a small $\eps'>0$ such that, for large enough $N,$ if $\| \fS(\tau) - \fe \|_{1} < \eps',$ then \begin{align*}
        \bE^{GC}[d(\sigma', \tau') | (\sigma, \tau)] \leq d(\sigma, \tau) \Big(1 - {\alpha \over N}\Big), 
    \end{align*}
    uniformly in $\sigma \in \Sigma.$
\end{prop}
\begin{proof}
    Write $\fz = \fS(\sigma)$ and $\fw = \fS(\tau)$ in $\cPm.$

    \textbf{Case 1.} Assume $\| \fz - \fe \|_{1} \geq \eps + \eps'$ and $\| \fw - \fe \|_{1} < \eps',$ where $\eps$ and $\eps'$ will be chosen sufficiently small. By Corollary \ref{D bound}, we have \begin{align*}
        \frac{D^{g}(\fz, \fe)}{\| \fz - \fe \|_{1}} \leq 1 - \gamma.
    \end{align*} By Remark \ref{monotone approx}, for sufficiently small $\eps>0,$ there is a coordinate-wise monotone path \begin{align*}
        \pi': \fe = \fz'_{0}, \ \fz'_{1}, \ \cdots, \ \fz'_{r} = \fz
    \end{align*} with $\| \fz'_{s-1} - \fz'_{s}\|_{1} \in (\eps, 2\eps]$ for each $s,$ and \begin{align*}
        D^{g}(\fz, \fe) = \sum_{i=1}^{m} \sum_{j=1}^{q} \sum_{s=1}^{r} | \inner{ \nabla g^{ij}_{\beta}(\fz'_{s-1})}{\fz'_{s} - \fz'_{s-1}}| + O(\eps^{2}).
    \end{align*}
    Since $\| \fw - \fe \|_{1}<\eps',$ we can adjust this path to a new monotone path \begin{align*}
        \pi: \fw = \fz_{0}, \ \fz_{1}, \ \cdots, \ \fz_{r}=\fz
    \end{align*} with $\| \fz_{s} - \fz'_{s}\|_{1} < \eps'$ for each $s.$ This path can be constructed by setting, for each $i \in [m], \ j \in [q],$ and $2 \leq s \leq r,$ \begin{align*}
    \fz_{s}^{ij} - \fz_{s-1}^{ij} \in \begin{cases}
        (0, (\fz_{s}')^{ij} - (\fz'_{s-1})^{ij}) & \text{if } (\fz_{s}')^{ij} - (\fz'_{s-1})^{ij}>0,\\
        ( (\fz_{s}')^{ij} - (\fz'_{s-1})^{ij}, 0) & \text{if } (\fz_{s}')^{ij} - (\fz'_{s-1})^{ij}<0.\\
    \end{cases}
    \end{align*}
    Then for each $i, \ j, \ s$, \begin{align*}
        \nonumber&| \inner{ \nabla g^{ij}_{\beta}(\fz'_{s-1})}{\fz'_{s} - \fz'_{s-1}} - \inner{ \nabla g^{ij}_{\beta}(\fz_{s-1})}{\fz_{s} - \fz_{s-1}} | \\
        \nonumber &\leq  | \inner{ \nabla g^{ij}_{\beta}(\fz'_{s-1})}{\fz'_{s} - \fz'_{s-1}} - \inner{ \nabla g^{ij}_{\beta}(\fz'_{s-1})}{\fz_{s} - \fz_{s-1}} |\\
        \nonumber&\quad + | \inner{ \nabla g^{ij}_{\beta}(\fz'_{s-1})}{\fz_{s} - \fz_{s-1}} - \inner{ \nabla g^{ij}_{\beta}(\fz_{s-1})}{\fz_{s} - \fz_{s-1}} |\\
         &=  O(\eps') + O(\eps' (\eps+\eps')).
    \end{align*}
   Applying Lemma \ref{contraction under monotone} to this path gives \begin{align*}
        \Gamma &\leq {1 \over \| \fz - \fw \|_{1}} \Big[\sum_{s=1}^{r} \sum_{i=1}^{m} \sum_{j=1}^{q} | \inner{ \nabla g^{ij}_{\beta} (\fz_{s-1})}{\fz_{s-1} - \fz_{s}} + O(\eps^{2}) + O\Big({\eps \over N}\Big) + O(N^{-2}) \Big] \\
        &\leq {1 \over \| \fz - \fw \|_{1}} \Big[\sum_{s=1}^{r} \sum_{i=1}^{m} \sum_{j=1}^{q} | \inner{ \nabla g^{ij}_{\beta} (\fz'_{s-1})}{\fz'_{s-1} - \fz'_{s}} + O(\eps') + O(\eps'(\eps+\eps'))+ O(\eps^{2}) + O\Big({\eps \over N}\Big) + O(N^{-2}) \Big]\\
        &\leq {1 \over \| \fz - \fw \|_{1}} \Big[D^{g}(\fz, \fe) + O(\eps') + O(\eps'(\eps+\eps'))+ O(\eps^{2}) + O\Big({\eps \over N}\Big) + O(N^{-2}) \Big]\\
        &\leq {1 \over \| \fz - \fe \|_{1}} \Big[D^{g}(\fz, \fe) + O(\eps') + O(\eps'(\eps+\eps'))+ O(\eps^{2}) + O\Big({\eps \over N}\Big) + O(N^{-2}) \Big] \Big[ 1 + O\Big({\eps' \over \eps+\eps'}\Big) \Big].
    \end{align*}
    By choosing $\eps' < \eps^{2}$ and $N > \eps^{-1},$ we have \begin{align*}
        \Gamma \leq 1- \gamma + O(\eps), 
    \end{align*}
    where constant $\gamma$ is defined in Corollary \ref{D bound}.
    
    \medskip

     \textbf{Case 2.}  Assume $\| \fz - \fe \|_{1} < \eps + \eps'$ and $\| \fw - \fe \|_{1} < \eps'.$ From $\eps' < \eps^{2},$ we have \begin{align*}
         \| \fz - \fw\|_{1} \leq \ \| \fz - \fe \|_{1} + \| \fw - \fw \| _{1} = O(\eps).
     \end{align*}
     Under the same Taylor approximation argument as in \eqref{q approx}, we get \begin{align*}
         \Gamma \leq & \frac{1}{\| \fz - \fw \|_{1}} \Big[ \sum_{i=1}^{m} \sum_{j=1}^{q} |\inner{\nabla g^{ij}_{\beta}(\fz)}{\fz - \fw}| + O(\eps^{2}) + O\Big( {\eps \over N} \Big) + O(N^{-2}) \Big] \\
         \leq & 1- \gamma_{0} + O(\eps),
     \end{align*} where $\gamma_0$ is from Remark \ref{when z, w are close}. 
     \medskip

    In both cases, by taking $\eps$ small enough, we obtain \begin{align*}
        \Gamma \leq 1 - \min\left\{ {\gamma \over 2}, {\gamma_{0} \over 2} \right\},
    \end{align*}
    which completes the proof.
\end{proof}

\medskip

\subsection{$O(N \log N)$ mixing}
To apply Proposition \ref{main contraction}, we need to estimate the probability that $\| \fS(\tau) - \fe\|_{1} < \eps'$ according to the Gibbs measure $\mu_{N, \beta}.$ We will show that this event happens with probability $1-O(e^{-cN})$ by a large deviation argument.

\begin{lemma} \label{LDP in upperbound}
    For sufficiently small $\eps'>0,$ define \begin{align*}
        A_{\eps', N} =\{ \sigma \in \Sigma: \| \fS(\sigma) - \fe \|_{1} < \eps' \}. 
    \end{align*} If $\beta<\beta_{c}(q)/J,$ then there exists a constant $c = c(\eps', \beta)>0$ such that \begin{align*}
        \mu_{N, \beta} [ A_{\eps', N}^{\co}] \leq \exp(-cN),
    \end{align*} 
    for all sufficiently large $N.$
\end{lemma}
\begin{proof}
    From \eqref{limsup}, we have \begin{align*}
        \mu_{N, \beta}[A^{\co}_{\eps', N}] = \nu_{N, \beta}[ \| \fS - \fe\|_{1} \geq \eps'] \leq \exp(-C N \min \{I_{\beta}(\fz): \|\fz - \fe\|_{1} \geq \eps'\} ),
    \end{align*}
    for all $C \in (0, 1).$ By Proposition \ref{phase trans}, $\min \{I_{\beta}(\fz): \|\fz - \fe\|_{1} \geq \eps'\}>0.$ Therefore, \begin{align*}
        \mu_{N, \beta} [ A_{\eps', N}^{c}] \leq \exp(-cN),
    \end{align*}
    for some constant $c = {1 \over 2} \min \{I_{\beta}(\fz): \|\fz - \fe\|_{1} \geq \eps'\}. $
\end{proof}

\setcounter{thm}{0}
\begin{thm}[Restated Theorem \ref{thrm 1}]
    Assume $\beta<\beta_s(q)/J.$ Then, for any $\eps>0,$ we have \begin{align*}
        t_{\operatorname{MIX}} \leq \Big({1 \over \alpha} + \eps\Big) N \log N,
    \end{align*}
    for sufficiently large $N,$ where $\alpha = \alpha(\beta)$ is the constant in Proposition \ref{main contraction}.
\end{thm}
\begin{proof}
    Let $(\sigma_{t}, \tau_{t})_{t \geq 0}$ be the greedy coupling with initial configurations $\sigma_0 \in \Sigma$ and $\tau_{0} \simeq \mu_{N, \beta}$ in distribution. Since $\tau_{t} \simeq \mu_{N, \beta}$ for all $t$, we have \begin{align*}
        \| \bP_{\sigma_0}[\sigma_t \in \cdot] - \mu_{N, \beta} \|_{\TV} \leq & \bP^{GC}[\sigma_t \neq \tau_t] = \bP^{GC}[d(\sigma_t, \tau_t) \geq 1]\\
        \leq & \bE^{GC}[d(\sigma_t, \tau_t)]\\
        \leq & \bE^{GC} [ \bE^{GC}[d(\sigma_t, \tau_t) | \tau_{t-1} \in A_{\eps', N} ] ] \ \bP^{GC}[\tau_{t-1} \in A_{\eps', N}]\\
        &+\diam(\Sigma) \ \bP^{GC}[\tau_{t-1} \notin A_{\eps', N}]\\
        \leq & \exp({-\alpha/N}) \bE^{GC} [d(\sigma_{t-1}, \tau_{t-1})] + \diam(\Sigma) \exp(-cN),
    \end{align*}
    where $\diam(\Sigma) = \max_{\sigma, \tau \in \Sigma} d(\sigma, \tau)= Nm$ and $c$ is as in Lemma \ref{LDP in upperbound}. We iterate this process by induction, so that \begin{align} \label{estimate}
         \|  \bP_{\sigma_0}[\sigma_t \in \cdot] - \mu_{N, \beta} \|_{\TV} \leq Nm (\exp({-\alpha t /N}) + t \exp(-cN)).
    \end{align}
    Put $t = (\alpha^{-1} + \eps) N \log N.$ Then, \eqref{estimate} becomes \begin{align*}
        Nm (e^{-(1+\alpha \eps) \log N} + (\alpha^{-1} + \eps) N \log N  \ e^{-cN}) \leq {m \over N^{\alpha \eps}} + e^{-c' N} \leq {1 \over 4},
    \end{align*}
   for sufficiently large $N,$ where $c'= c'(\alpha, \eps)>0.$ Therefore, $t_{\mix} \leq (\alpha^{-1} + \eps) N \log N$ as claimed. 
\end{proof}

\section{Exponential mixing in the low-temperature regime $\beta> \beta_{s}(q)/J$} \label{low}
Since the metastability of the free energy implies exponential mixing \cite{markov, Glauber, yang}, it suffices to establish metastability in the low-temperature regime $\beta > \beta_s(q) / J.$ We will prove that there are at least $q$ local minimizers of the free energy $F_{\beta}$ defined in \eqref{freedef}. Inspired by \cite{landscapepotts, twocomp}, we show that the matrix \begin{align*}
   \begin{bmatrix}
\tilde{t}_2 \ \; \ \tilde{t}_2 \; \ \tilde{t}_2 \; \ \cdots \ \; \ 1-(q-1) \tilde{t}_{2}\\
\tilde{t}_2 \ \; \ \tilde{t}_2 \; \ \tilde{t}_2 \; \ \cdots \ \; \ 1-(q-1) \tilde{t}_{2}\\
\vdots\\
\tilde{t}_2\ \; \ \tilde{t}_2 \; \ \tilde{t}_2 \; \ \cdots \ \; \ 1-(q-1) \tilde{t}_{2}
\end{bmatrix} \in \cPm, \quad \text{ with } \tilde{t}_{2}(\beta)= \frac{1-x_{2}(\beta J)}{q-1},
\end{align*} where $x_{2}(\beta)$ is from Lemma \ref{solution study}, is a local minimizer. We begin by analyzing $$t_{2}(\beta) = \frac{1-x_{2}(\beta)}{q-1}$$ when $\beta>\beta_{s}(q).$ 
 \begin{lemma} \label{meta lemma}
    Let $\beta> \beta_{s}(q).$ Then, $t_{2}(\beta) \in (0, 1/q),$ and \begin{align*}
        2\beta < \min \left\{{1 \over t_{2}(\beta)}, {1 \over q t_{2}(\beta) (1-(q-1)t_{2}(\beta))}\right\}.
    \end{align*}
 \end{lemma}
 \begin{proof}
     Since $x_{2}(\beta) \in (1/q, 1)$, it follows that $t_{2}(\beta) \in (0, 1/q).$ Moreover, \begin{align*}
         \frac{t_{2}(\beta)}{1-(q-1) t_{2}(\beta)} = {1 \over q-1} \ \frac{1-x_2(\beta)}{x_2(\beta)} =\exp({2\beta(qt_2(\beta) -1)}).
     \end{align*}
     Therefore, $t_{2}(\beta) \in (0, 1/q)$ is a solution of the equation \begin{align} \label{equation}
         2\beta =\xi(t):= {1 \over 1-qt} \log \Big(\frac{1-(q-1)t}{t}\Big).
     \end{align}
     By similar arguments, $t_{1}(\beta)=(1-x_{1}(\beta))/(q-1)$ is the other solution solution of \eqref{equation} with $t_{1}(\beta)>t_{2}(\beta).$ Computing \begin{align*}
         \xi'(t)=\frac{q}{(1-qt)^{2}} \Big[ \log \Big(\frac{1-(q-1)t}{t}\Big) + \frac{qt-1}{qt(1-(q-1)t)} \Big],\quad
         \xi''(t)=\frac{1- 2(q-1)t}{(1-qt)(1-(q-1)t)^{2} t^{2}},
     \end{align*}
     one sees that $\xi$ has a unique minimizer on $(0, 1/q).$ Thus, $2\beta = \xi(t)$ has exactly two solutions $t_{1}(\beta)>t_{2}(\beta)$ with $\xi'(t_2(\beta))<0<\xi'(t_1(\beta)).$

     To show $2\beta<{1/ t_{2}(\beta)},$ note \begin{align*}
         2\beta<{1 \over t_{2}(\beta)} \iff {1 \over t_{2}(\beta)} > {1 \over 1-qt_{2}(\beta)} \log \Big(\frac{1-(q-1)t_{2}(\beta)}{t_{2}(\beta)}\Big) \iff \frac{1-qt_{2}(\beta)}{t_{2}(\beta)}>\log \Big(1 + \frac{1-qt_{2}(\beta)}{t_{2}(\beta)}\Big),
     \end{align*}
     since $(1-qt_{2}(\beta))/t_{2}(\beta) >0.$

    In addition, since $\xi'(t_{2}(\beta))<0,$ we obtain \begin{align*}
         &\log\Big(\frac{1-(q-1)t_{2}(\beta)}{t_{2}(\beta)}\Big) < \frac{1-qt_{2}(\beta)}{qt_{2}(\beta) (1-(q-1)t_{2}(\beta))}\\ &\iff 2\beta =  {1 \over 1-qt_{2}(\beta)} \log\Big(\frac{1-(q-1)t_{2}(\beta)}{t_{2}(\beta)}\Big) < \frac{1}{qt_{2}(\beta) (1-(q-1)t_{2}(\beta))}.
     \end{align*}
 \end{proof}
\begin{prop}
    If $\beta> \beta_{s}(q)/J,$ then $F_{\beta}$ has at least $q$ local minimizers in $\cPm.$
\end{prop}
\begin{proof}
    Write the free energy as \begin{align*}
        F_{\beta}(\fS) =  -\beta \Tr(\fz^{\intercal} \fJ \fz) + \sum_{i=1}^{m} \sum_{j=1}^{q} \fz^{ij} \log(\fz^{ij}), \quad \fS \in \cPm,
    \end{align*}
    and regard it as a function in $\bR^{m \times (q-1)}$ by considering \begin{align*}
        \fS^{kq} = 1 - \fS^{k1} - \fS^{k2} - \cdots - \fS^{k(q-1)}, \ \text{ for each } k \in [m].
    \end{align*}
    Accordingly, we define \begin{align*}
        \fT = \begin{bmatrix}
	\tilde{t}_2 \ \; \ \tilde{t}_2 \; \ \tilde{t}_2 \; \ \cdots \ \; \ \tilde{t_{2}}\\
\tilde{t}_2 \ \; \ \tilde{t}_2 \; \ \tilde{t}_2 \; \ \cdots \ \; \ \tilde{t}_{2}\\
\vdots\\
\tilde{t}_2\ \; \ \tilde{t}_2 \; \ \tilde{t}_2 \; \ \cdots \ \; \ t_{2}
\end{bmatrix} \in \bR^{m \times (q-1)}, \quad \fT^{iq} = 1 - (q-1) \tilde{t}_{2}(\beta), \ \forall i \in [m], 
    \end{align*}
    where $\tilde{t}_{2} = t_{2}(\beta J) = (1 - x_{2}(\beta J))/(q-1).$ One checks that for each $k \in [m]$ and $l \in [q-1],$ \begin{align*}
        \frac{\partial F_{\beta}}{\partial \fS^{kl}}=- 2\beta J (\fS^{kl} - \fS^{kq}) - 2\beta J \lambda \sum_{u \neq k}^{m} (\fS^{ul} - \fS^{uq}) + (1 + (m-1)\lambda) (\log \fS^{kl} - \log \fS^{kq}).
    \end{align*}
    Hence $\fT$ is a critical point, because \begin{align*}
        \frac{\partial F_{\beta}}{\partial \fS^{kl}}(\fT)= 2\beta J (1 + (m-1) \lambda) \Big[ (1-q\tilde{t}_2) - {1 \over 2\beta J} \log \Big( \frac{1-(q-1)\tilde{t}_2}{\tilde{t}_2} \Big)\Big] = 0.
    \end{align*}
    We now examine the second derivatives of $F_{\beta}.$ For distinct $k, r \in [m]$ and distinct $l, p \in [q-1],$ we have \begin{align*}
        &\frac{\partial^{2} F_{\beta}}{\partial (\fS^{kl})^{2}} (\fT) = 2\beta J( a + b), \quad \frac{\partial^{2} F_{\beta}}{\partial \fS^{kl} \partial \fS^{kp}}(\fT) = 2\beta J b, \\
        &\frac{\partial^{2} F_{\beta}}{\partial \fS^{kl} \partial \fS^{rl}} (\fT)= -4 \beta J \lambda, \quad \frac{\partial^{2} F_{\beta}}{\partial \fS^{kl} \partial \fS^{kp}} (\fT)= - 2 \beta J \lambda,      
    \end{align*}
    where \begin{align*}
        a = -1 + \frac{1 + (m-1)\lambda}{2\beta J \tilde{t}_2}, \quad b = -1 + \frac{1 + (m-1)\lambda}{2\beta J( 1 - (q-1) \tilde{t}_2)}.
    \end{align*}
    Identifying $\bR^{m \times (q-1)} \simeq \bR^{m(q-1)},$ the Hessian has block form \begin{align*}
        \nabla^{2} F_{\beta} (\fT) = 2 \beta J\begin{bmatrix}
            \mathbf{M} \ \; \ \mathbf{N} \; \ \mathbf{N} \; \ \cdots \ \; \ \mathbf{N}\\
\mathbf{N} \ \; \ \mathbf{M} \; \ \mathbf{N} \; \ \cdots \ \; \ \mathbf{N}\\
\vdots\\
\mathbf{N} \ \; \ \mathbf{N} \; \ \mathbf{N} \; \ \cdots \ \; \ \mathbf{M}\\
        \end{bmatrix} \in \bR^{m(q-1) \times m(q-1)}
    \end{align*}
    with on-diagonal block $\mathbf{M} \in \bR^{(q-1) \times (q-1)}$ and off-diagonal block $\mathbf{N} \in \bR^{(q-1) \times (q-1)}$ such that \begin{align*}
        \mathbf{M} = a \fI_{q-1} + b \mathbf{1}_{q-1}, \quad \mathbf{N} = -\lambda \fI_{q-1} -\lambda \mathbf{1}_{q-1}, 
    \end{align*} where $\mathbf{I}_{q-1}$ is $(q-1) \times (q-1)$ identity matrix and $\mathbf{1}_{q-1}$ is $(q-1) \times (q-1)$ matrix with all entries one. As a result, $\nabla^{2} F_{\beta}(\fT)$ has the eigenvalues \begin{enumerate}
        \item $a + (q-1)b + q\lambda,$ with multiplicity $m-1,$
        \item $a + (q-1)b - q(m-1)\lambda,$ with multiplicity $1,$
        \item $a +\lambda,$ with multiplicity $(q-2)(m-1),$
        \item $a - (m-1) \lambda,$ with multiplicity $q-2.$
    \end{enumerate}
    For positive-definiteness, it is enough to show $a - (m-1) \lambda >0$ and $a+b(q-1) - q(m-1)\lambda>0.$ Indeed, \begin{align}
        a > (m-1) \lambda  \iff 2\beta J <{1 \over \tilde{t}_{2}}, \quad a + (q-1)b > q(m-1) \lambda \iff 2\beta J <{1 \over q \tilde{t}_{2}(1 - (q-1)\tilde{t}_{2})}, 
    \end{align}
   which follows from Lemma \ref{meta lemma}. So $\fT$ is a local minimizer. By permuting the columns of $\fT,$ we obtain $q$ different local minimizers.
\end{proof}

\bigskip

\section*{Additional Information}
\noindent Data availability: No datasets were generated or analyzed during the current study.\\
Competing interests: The author declares no competing financial interests.

\bigskip

\bibliographystyle{alpha}
\bibliography{references.bib}


\end{document}